\pdfoutput=1
\documentclass[12pt,number]{elsarticle}

\usepackage{amssymb}
\usepackage{latexsym}
\usepackage{color}
\usepackage{graphicx}
\newcommand{\trunc}[2]{\mbox{trunc}_{#1}\left(#2\right)}

\begin{document}

\begin{frontmatter}



\title{Matrix approach to discrete fractional calculus~II:
       partial fractional differential equations}

\author[tuke]{Igor Podlubny\corref{cor1}} 
\ead{igor.podlubny@tuke.sk} 
\author[kharkov]{Aleksei Chechkin} 
\ead{achechkin@kipt.kharkov.ua} 
\author[tuke]{Tomas Skovranek} 
\ead{tomas.skovranek@tuke.sk} 
\author[usu]{YangQuan Chen} 
\ead{yqchen@ieee.org} 
\author[unex]{Blas M. Vinagre Jara} 
\ead{bvinagre@unex.es}

\cortext[cor1]{Corresponding author}
       
\address[tuke]{BERG Faculty, Technical University of Kosice, Slovak Republic}
\address[kharkov]{Institute for Theoretical Physics NSC KIPT, Kharkov, Ukraine}
\address[usu]{Department of Electrical and Computer Engineering, Utah State University, USA}
\address[unex]{Industrial Engineering School, University of Extremadura, Badajoz, Spain}

\begin{abstract}
A new method that enables easy and convenient discretization of partial differential equations with 
derivatives of arbitrary real order (so-called fractional derivatives) and delays
is presented and illustrated on numerical solution
of various types of fractional diffusion equation. The suggested method is the development
of Podlubny's matrix approach (Fractional Calculus and Applied Analysis, vol. 3, no. 4, 2000, 359--386).
Four examples of numerical solution of fractional diffusion equation with various combinations of time/space fractional derivatives (integer/integer, fractional/integer, integer/fractional, and fractional/fractional) with respect to time and to the spatial variable are provided in order to illustrate how simple and general is the suggested approach.  The fifth example illustrates 
that the method can be equally simply used for fractional differential equations with delays. 
A set of MATLAB routines for the implementation of the method as well as sample code 
used to solve the examples have been developed.
\end{abstract}

\end{frontmatter}

% main text
\section{Introduction}

%%% Physical motivation

Recently, kinetic equations of the diffusion, diffusion--advection, and \linebreak Fokker--Planck type with partial fractional derivatives were recognized as a useful approach for the description of transport dynamics in complex systems whose temporal evolution deviates from the standard laws, that is, from exponential Debye or Gaussian laws, and from fast decaying correlations. Examples include systems exhibiting Hamiltonian chaos, disordered medium, plasma and fluid turbulence, underground water pollution, dynamics of protein molecules, motions under the influence of optical tweezers, reactions in complex systems, and more (see reviews on fractional kinetics~\cite{CheGoMeKla-ACP06,MeKla-PhysRep00,SoKlaBlu-PhysTo02,Zaslav-PhysRep02,MeKla-JPA03}, the recent multi-author book \cite{Klagers-2008},  and references therein). These fractional equations are derived asymptotically from basic random walk models, the generalized master and Langevin equations. The advantage of the fractional models lies in the straightforward way of including external force terms and of calculating boundary value problems. Also, the consideration of transport in the phase space spanned by both position and velocity coordinates is possible within the fractional approach. However, because of complicated integro-differential structure of fractional kinetic equations the analytical solutions are presently known only in a very few relatively simple cases. Therefore, the development of numerical methods is of current importance.

%%%%%%%%%%

%%% Types of ``fractionalization'' 

Let us recall briefly how the kinetic equations with integer partial derivatives can be ``fractionalized''. There are two generic types of fractionalization, which can be explained by taking as an example the parabolic diffusion equation for the particles density $u(x,t)$ in a one-dimensional space,

\begin{equation}\label{eq:classic-diffusion}
\frac{\partial u}{\partial t} =
\chi 
\frac{\partial^2 u}{\partial x^2},
\quad (\;t > 0, \;\;  a < x < b \;)
\end{equation}

\noindent
where constant $\chi$  is diffusion coefficient. The first type of fractionalization leads to a \emph{time fractional diffusion equation} by means of replacing the first order time derivative by afractional derivative of  order $\alpha$  less than 1,

\begin{equation}\label{eq:time-frac-diffusion}
_{0}^{C}\!D_{t}^{\alpha} u =
\chi 
\frac{\partial^2 u}{\partial x^2},
\quad (\;t > 0, \;\;  a < x < b \;)
\end{equation}

Here,  $_{0}^{C}\!D_{t}^{\alpha} $ is the Caputo fractional derivative~\cite{Caputo-1969}, which is defined as

\begin{equation}\label{eq:Caputo-definition-left}
    _{a}^{C}\!D_{x}^{\mu} \phi(x) =
    \frac{1}{\Gamma (m -\mu)}
    \int\limits_{a}^{x}
    \frac{\phi^{(m)}(\xi) d\xi}
         {(x - \xi)^{\mu - m + 1}},
\qquad
( m -1 < \mu \leq m )
\end{equation}

Taking $\alpha = 1$ in (\ref{eq:time-frac-diffusion}) gives the classical diffusion equation (\ref{eq:classic-diffusion}).

Other two forms of a time fractional diffusion equation that appears in the literature use the Riemann-Liouville fractional derivative instead of the Caputo one~\cite{MeKla-PhysRep00}. 
Although recently, in addition to a geometric and physical interpretation of fractional integration and fractional differentiation \cite{Podlubny-geomint-2002}, a physical interpretation for the initial conditions in terms of the Riemann-Liouville fractional derivatives of the unknown function has been suggested \cite{Heymans:2006fy},
the use of Caputo derivative in physical problems is perhaps more convenient since it allows using initial conditions expressed in terms of values of the unknown function and  its integer-order derivatives~\cite{Podlubny-FDE-book}. However, all three forms of ``time-fractionalization'' are equivalent if zero initial conditions are posed. In what follows we use the form with the Caputo derivative, equation (\ref{eq:time-frac-diffusion}), since some of the illustrating examples use non-zero initial conditions.

In the second type of fractionalization, the second order spatial derivative is replaced by the fractional derivative of the order $\beta$ between 1 and 2, thus leading to spatial fractional diffusion equation,

\begin{equation}\label{eq:intro-1-5}
\frac{\partial u}{\partial t} =
\chi 
\frac{\partial^\beta u}{\partial |x|^\beta},
\quad (\;t > 0, \;\;  a < x < b \;)
\end{equation}

\noindent
where $\partial^{\beta}/\partial |x|^{\beta}$ (we adopt here the notation introduced in~\cite{SZ97}) is a partial (with respect to spatial variable) symmetric Riesz derivative, which is defined as a half-sum of the left- and right-sided Riemann-Liouville derivatives \cite{Podlubny-FDE-book,Podlubny-MFC-2000}:

\begin{equation}\label{eq:Riesz-derivative-definition}
\frac{d^\beta \phi(x)}{d|x|^\beta} =
D_{R}^{\beta} \phi(x) = 
\frac{1}{2}
\Bigl(  \, \,
_{a}D_{x}^{\beta} \phi(x) + \, 
_{x}D_{b}^{\beta} \phi(x)
\Bigr),
\end{equation}
\noindent
where the left- and right-sided Riemann-Liouville derivatives are defined by 

\begin{equation}\label{eq:RL-definition-left}
    _{a}D_{x}^{\mu} \phi(x) =
    \frac{1}{\Gamma (m -\mu)}
    \left(
        	\frac{d}{dx}
    \right)^m
    \int\limits_{a}^{x}
    \frac{\phi (\xi) d\xi}
         {(x - \xi )^{\mu - m + 1}},
\quad
( m -1 < \mu \leq m ),
\end{equation}

\begin{equation}\label{eq:RL-definition-right}
    _{x}D_{b}^{\mu} \phi(x) =
    \frac{1}{\Gamma (m -\mu)}
    \left(
       -
    	\frac{d}{dx}
    \right)^m
    \int\limits_{x}^{b}
    \frac{\phi (\xi) d\xi}
         {(\xi - x)^{\mu - m + 1}},
\quad
( m -1 < \mu \leq m ),
\end{equation}

For $\beta = 2$ the equation (\ref{eq:intro-1-5}) becomes the classical diffusion equation (\ref{eq:classic-diffusion}).

Other forms of asymmetric space fractional generalizations use the left-side Riemann - Liouville derivative instead of the symmetric Riesz derivative~\cite{Diego-PRL03,MeerTad04}, or asymmetric derivative with different asymmetry factors at the left- and right-side derivatives~\cite{Diego-PRL05,DiegoPhysPlas,Meerschaert2006}. In terms of random walk schemes, the symmetric derivative corresponds to a symmetric jump probability distribution of a diffusing particle, whereas any asymmetry in space derivative accounts for inherent force-free preferable direction of jumps which may occur, e.g., in heterogeneous porous media or magnetically confined fusion plasmas. In our paper we restrict ourselves to symmetric case, equation (\ref{eq:intro-1-5}).

Of course, there are different generalizations of time and space fractional diffusion equations, including: multidimensional fractional diffusion and kinetic equations~\cite{ChePlas,Friedrich03,Meer1999}, both time and space fractional generalizations~\cite{MaiLu01}, different regular forces in space and time fractional Fokker-Planck equations~\cite{MeBarKla99,CheChemPhys,Added-11,Added-12,Added-13,Added-14}, variable transport coefficients~\cite{Meer06}, equations with fractional derivatives of distributed and variable orders
~\cite{PhysRevE02,JPA05,Added-21,Added-22} etc. The realm of fractional kinetics is growing, and therefore it is desirable to have at hand a method for numerical solution which would be relatively simple and at the same time general enough to deal effectively with different forms of fractional kinetic equations. However, while different numerical tools for ordinary fractional equations exist and a basic framework of their numerical solution is already established, relatively few numerical methods exist to solve fractional equations with partial derivatives, and the development of effective numerical schemes is now on the agenda. We recall briefly the different approaches used in the literature.

The numerical methods differ essentially in the way in which normal and fractional derivatives are discretized.
In \cite{Lynch-et-al-2003} to solve diffusion-reaction equation with the left Riemann-Liouville derivative between 1 and 2, the L2 discretization method was used taken 
from~\cite{Oldham-Spanier}, together with its modification, L2C (both L2 and L2C methods are based on numerical approximation of a fractional integral that appears in the definition of the Riemann-Liouville fractional derivative). 
It was shown that the former is the most accurate for orders larger than 1.5, whereas the latter is the most accurate for orders less than~1.5. For the first order time derivative, the explicit forward Euler formula and semi-implicit scheme were used. 

Langlands and Henry~\cite{Langlands-2005} used L1 scheme from~\cite{Oldham-Spanier} to discretize the Riemann-Liouville fractional time derivative of order between 1 and 2.

Yuste~\cite{Yuste2006} considered a Gr\"unwald-Letnikov approximation for the Riemann-Liouville time derivative and used a weighted average for the second-order space derivative.

Scherer et al. \cite{Scherer:2008kk} introduced very recently a modification of the 
Gr\"unwald-Letnikov approximation for the case of the Caputo derivative of a function
which is not zero in the starting point of the considered time interval, and applied 
that approximation for the numerical solution of fractional diffusion equations
with the Caputo time derivative and non-zero initial conditions.

To solve the one-dimensional space fractional advection-dispersion equation with left-side Riemann-Liouville derivative and variable coefficients the shifted 
Gr\"un\-wald-Letnikov approximation was proposed by Meerschaert and Tadjeran~\cite{MeerTad04}. For two-sided space-fractional partial differential equations the shifted Gr\"un\-wald-Letnikov formula was proposed and discussed in~\cite{MeerTad06}. The fractional Crank-Nicholson method based on the shifted formula was elaborated, giving temporally and spatially second-order numerical estimates~\cite{TadMeer06}. The generalizations of the shifted formula and of the fractional Crank-Nicholson method in the two-dimensional case were discussed in~\cite{MeerSchef06} and ~\cite{TadMeer07}, respectively.

Another method to solve the space-fractional Fokker-Planck equation with constant coefficient on the fractional derivative term was pursued by Liu et al. ~\cite{Liu2004}. They transform the partial differential equation into a system of ordinary differential equations, which is solved by a method of lines.

Ervin and Roop ~\cite{Ervin2005,Ervin2006} presented a theoretical framework for the Galerkin finite element approximation to the steady state fractional advection-diffusion equation, and extended this approach to multidimensional partial differential equations with constant coefficients on the fractional derivative terms.

Valko and Abate~\cite{Valko2005}  solved the time-fractional diffusion equation on a semi-infinite domain by numerical inversion of the two-dimensional Laplace transform. To solve the time-fractional diffusion equation in a bounded domain, Lin and Xu ~\cite{Lin07} proposed the method based on a finite difference scheme in time and Legendre spectral method in space.

Liang and Chen \cite{Liang:2006nx} used a combination of symbolic computations and numerical inversion of the Laplace transform for solving a time-fractional diffusion-wave equation with the time derivative of order between 1 and 2. 

We also mention that in order to approximate shifted Caputo time derivative appearing in hydrodynamic equations for heterogeneous porous media the modification of Yuan and Agrawal's method ~\cite{Yuan2002} was used to transform a fractional derivative into an infinite integral over auxiliary internal variables ~\cite{Lu2005}.

%%%% 

Another approach for the solution of fractional kinetic equations employs the methods of Monte Carlo type (random walk based methods).
A set of random walk schemes applied to fractional diffusion equations based on the Gr\"unwald-Letnikov approximation was developed in the papers by Gorenflo, Mainardi and co-workers. They were applied to solve (i) symmetric space-fractional diffusion equation ~\cite{Gorenflo1999,Gorenflo2001}; (ii) asymmetric space-fractional diffusion equation in the L\`evy--Feller form ~\cite{Gorenflo1998}; (iii) time-fractional diffusion equation with Caputo derivative~\cite{GoMa02}; (iv) time-space fractional diffusion equation ~\cite{GoMaMo02,GoMaMo02a}. Chechkin et al. ~\cite{CheFCAA03} generalized the approach on a double-order time fractional diffusion equation. Gorenflo and Abdel-Rehim ~\cite{GoAb05} proposed discrete approximations to time-fractional diffusion process with non-homogeneous drift towards the origin by generalization of Ehrenfest's urn model. The L\`evy--Feller diffusion-advection process with a constant drift was approximated by random walk and finite difference method by Liu et al. ~\cite{Liu2007}. The random walk particle tracking approach to solve one-dimensional space-fractional diffusion-advection equation with space dependent coefficients was employed by Meerschaert and co-authors ~\cite{Meer06}. The method based on numerical solution of a coupled stochastic differential equations driven by L\`evy symmetric stable processes was proposed in~\cite{Stanescu2005} to solve a non-linear evolution problem involving the fractional Laplacian operator.

All aforementioned works indicate that numerical solution of partial fractional differential equations plays an important and increasing role in the applications of the methods and models of non-integer order.

In the present paper we propose a general approach to the numerical solution of  partial fractional differential equations, which is based on the matrix form representation of discretized fractional operators introduced in \cite{Podlubny-MFC-2000}. This approach unifies the numerical differentiation of arbitrary (including integer) order and the $n$-fold integration, using the so-called triangular matrices. Applied to numerical solution of differential equations, it also unifies the solution of integer- and fractional-order partial differential equations.
The suggested approach leads to significant simplification of the numerical solution of partial differential equations, and it is general enough to deal with different types of partial fractional differential equations, even with delays.

\section{The idea of the suggested method}

The method that we suggest is based on triangular strip matrix approach ~\cite{Podlubny-MFC-2000} to discretization of operators of differentiation and integration of arbitrary real order.

\begin{figure}
\begin{center}
\includegraphics[width=0.6\textwidth]{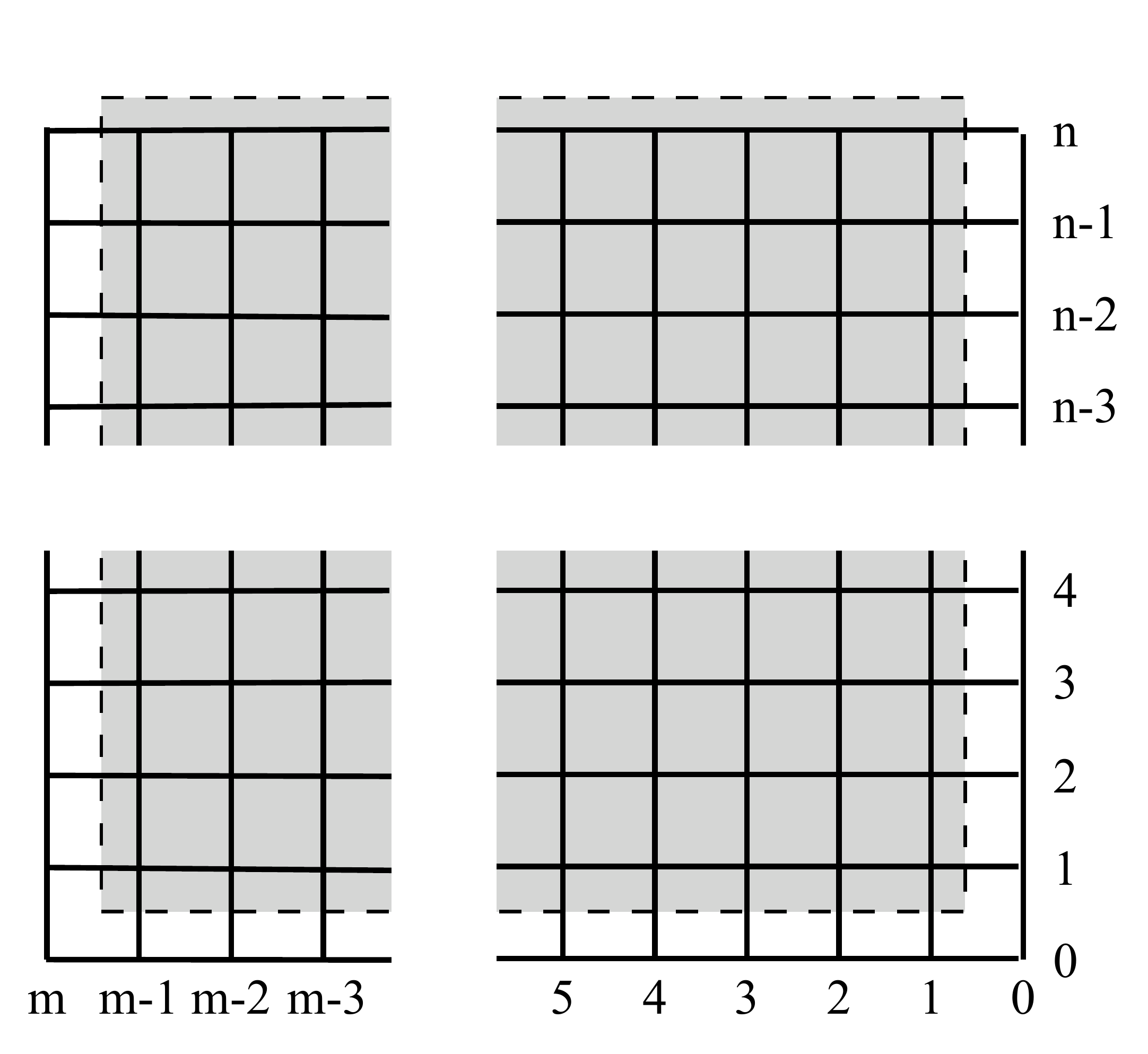}
\caption{Nodes and their right-to-left, and bottom-to-top numbering.}
\label{fig:nodes-net}
\end{center}
\end{figure}

In contrast with generally used numerical methods, where the solution is obtained step-by-step by moving from the previous time layer to the next  one, let us consider the whole time interval 
of interest at once. 
This allows us to create a net of discretization nodes. In the simplest case of one spatial dimension this step gives a 2D net of nodes. An example of such discretization is shown in Fig.~\ref{fig:nodes-net}. 
The values of the unknown function in inner nodes (shaded area in Fig.~\ref{fig:nodes-net}) are to be found. The values at the boundaries are known:  they are used later in constructing the system of algebraic equations.

The system of algebraic equations is obtained by approximating the equation in all inner nodes simultaneously (this gives the left-hand side of the resulting system of algebraic equations) and then utilizing the initial and boundary conditions (the values of which appear in the right-hand side of the resulting system). 

The discretization nodes in Fig.~\ref{fig:nodes-net} are numbered from right to left in each time level, and the time levels are numbered from bottom to top. We use such numbering in this article for the clarity of presentation of our approach, although standard numberings work equally well.

In the following sections we recall the basic tools that are necessary for the method: the triangular strip matrices, the Kronecker product, the eliminators, and the shifters. Then we show how they are used for approximating partial derivatives of arbitrary real order and the equation, and how the resulting system of algebraic equations appears.

\section{Triangular strip matrices}

In this paper we use matrices of a specific structure,
which are called {\em triangular strip matrices\/}
\cite{Podlubny-MFC-2000,Suprunenko-Tyshkevich-66}, 
and which have been also mentioned
in \cite{Bulgakov-54,Gantmakher-Matrix-Theory}.
We will need lower triangular strip matrices,
\begin{equation}\label{eq:General-lower-triangular}
    L_N=
    \left[
    \begin{array}{cccccc}
    \omega_0   & 0  & 0  & 0  & \cdots & 0 \\
    \omega_1   & \omega_0   & 0  & 0  & \cdots & 0 \\
    \omega_2   & \omega_1  & \omega_0 & 0  & \cdots & 0 \\
    \ddots & \ddots & \ddots & \ddots & \cdots & \cdots \\
    \omega_{N-1}    & \ddots  & \omega_2   & \omega_1 & \omega_0  & 0 \\
    \omega_N  & \omega_{N-1}   & \ddots  & \omega_2  & \omega_1  & \omega_0  \\
    \end{array}
    \right],
\end{equation}

\noindent
and upper triangular strip matrices,
\begin{equation}\label{eq:General-upper-triangular}
    U_N=
    \left[
    \begin{array}{cccccc}
    \omega_0   & \omega_1  & \omega_2  & \ddots  & \omega_{N-1} & \omega_N \\
    0   & \omega_0   & \omega_1   & \ddots  & \ddots & \omega_{N-1}  \\
    0   & 0 & \omega_0   & \ddots   & \omega_2    & \ddots  \\
    0   & 0 & 0          & \ddots   & \omega_1    & \omega_2 \\
    \cdots  & \cdots     & \cdots   & \cdots      & \omega_0  & \omega_1 \\
    0   & 0 & 0    & \cdots & 0    & \omega_0  \\
    \end{array}
    \right],
\end{equation}

A lower (upper) triangular strip matrix is completely described by its first column (row). 
Therefore, if we define the truncation operation, $\trunc{N}{\cdot}$, which truncates
(in a general case) the power series $\varrho(z)$,
\begin{equation}\label{eq:Generating-function-general}
    \varrho(z)=\sum_{k=0}^{\infty} \omega_k z^k
\end{equation}
to the polynomial $\varrho_N(z)$,
\begin{equation}
    \trunc{N}{\varrho(z)}
    \stackrel{\mbox{\scriptsize def}}{=}
    \sum_{k=0}^{N} \omega_k z^k=
    \varrho_N(z),
\end{equation}
then we can consider the function $\varrho(z)$ as a generating series
for the set of lower (or upper) triangular matrices $L_N$ (or $U_N$),
$N=1, 2, \ldots $

It was shown in  ~\cite{Podlubny-MFC-2000}  that operations with triangular strip matrices, 
such as addition, subtraction, multiplication, and inversion, can be expressed in 
the form of operations with their generating series (\ref{eq:Generating-function-general}).

Among properties of triangular strip matrices it should be noticed that
 if matrices $C$ and $D$ are both lower (upper) triangular
strip matrices, then they commute:
\begin{equation}
    CD = D \, C.
\end{equation}

%%%%%%%  KRONECKER PRODUCT %%%%%%%%%%%%%%%%%%%%%%

\section{Kronecker matrix product}

The Kronecker product $A \otimes B$ of the  $n\times m$ matrix $A$ and the $p \times q$ matrix $B$
\begin{equation}
A =
\left[
	\begin{array}{cccc}
			a_{11} & a_{12} & \ldots & a_{1m} \\
			a_{21} & a_{22} & \ldots & a_{2m} \\
			\vdots & \vdots & \ddots & \vdots \\
			a_{n1} & a_{n2} & \ldots & a_{nm}
	\end{array}
\right], 
\qquad 
B =
\left[
	\begin{array}{cccc}
			b_{11} & b_{12} & \ldots & b_{1q} \\
			b_{21} & b_{22} & \ldots & b_{2q} \\
			\vdots & \vdots & \ddots & \vdots \\
			b_{p1} & b_{p2} & \ldots & b_{pq}
	\end{array}
\right], 
\end{equation}
is the $np \times mq$ matrix having the following block structure:
\begin{equation}
A \otimes B=
\left[
	\begin{array}{cccc}
			a_{11}B & a_{12}B & \ldots & a_{1m}B \\
			a_{21}B & a_{22}B & \ldots & a_{2m}B \\
			\vdots & \vdots & \ddots & \vdots \\
			a_{n1}B & a_{n2}B & \ldots & a_{nm}B
	\end{array}
\right]. 
\end{equation}

For example, if 
\begin{equation}
A = 
\left[\begin{array}{cc}1 & 2 \\  0 & -3\end{array}\right], 
\qquad
B = 
\left[\begin{array}{ccc}1 & 2 & 3 \\4 & 5 & 6\end{array}\right],
\end{equation}
then
\begin{equation}
A \otimes B = 
\left[\begin{array}{rrrrrr}1 & 2 & 3 & 2 & 4 & 6 \\4 & 5 & 6 & 8 & 10 & 12 \\0 & 0 & 0 & -3 & -6 & -9 \\0 & 0 & 0 & -12 & -15 & -18\end{array}\right].
\end{equation}

Among many known interesting properties of the Kronecker product we would like to recall
those that are important for the subsequent sections. Namely~\cite{Loan:2000},
\begin{itemize}
\item 
if $A$ and $B$ are band matrices, then $A \otimes B$ is also a band matrix,
\item
if $A$ and $B$ are lower (upper) triangular, then $A \otimes B$ is also  lower (upper) triangular.
\end{itemize}

We will also need two specific Kronecker products, namely the products 
$E_{n} \otimes A$ and $A \otimes E_{m}$, where $E_{n}$ is an $n \times n$ identity matrix.
For example, if $A$ is a $2 \times 3$ matrix 

\begin{equation}
A = 
 \left[ \begin {array}{ccc} {\it a_{11}}&{\it a_{12}}&{\it a_{13}}\\\noalign{\medskip}{\it a_{21}}&{\it a_{22}}&{\it a_{23}}\end {array} \right] 
\end{equation}

\noindent
then 

\begin{equation}
E_{2} \otimes A = 
\left[ \begin {array}{cccccc} {\it a_{11}}&{\it a_{12}}&{\it a_{13}}&0&0&0\\\noalign{\medskip}{\it a_{21}}&{\it a_{22}}&{\it a_{23}}&0&0&0\\\noalign{\medskip}0&0&0&{\it a_{11}}&{\it a_{12}}&{\it a_{13}}\\\noalign{\medskip}0&0&0&{\it a_{21}}&{\it a_{22}}&{\it a_{23}}\end {array} \right] 
\end{equation}

\begin{equation}
A \otimes E_{3} = 
 \left[ \begin {array}{ccccccccc} {\it a_{11}}&0&0&{\it a_{12}}&0&0&{\it a_{13}}&0&0\\\noalign{\medskip}0&{\it a_{11}}&0&0&{\it a_{12}}&0&0&{\it a_{13}}&0\\\noalign{\medskip}0&0&{\it a_{11}}&0&0&{\it a_{12}}&0&0&{\it a_{13}}\\\noalign{\medskip}{\it a_{21}}&0&0&{\it a_{22}}&0&0&{\it a_{23}}&0&0\\\noalign{\medskip}0&{\it a_{21}}&0&0&{\it a_{22}}&0&0&{\it a_{23}}&0\\\noalign{\medskip}0&0&{\it a_{21}}&0&0&{\it a_{22}}&0&0&{\it a_{23}}\end {array} \right] 
\end{equation}

This illustrates that left multiplication  of $A_{n \times m}$ by $E_{n}$ 
creates an $n \times n$ block diagonal matrix by repeating 
the matrix $A$ on the diagonal, 
and that right multiplication of $A_{n \times m}$ by $E_{m}$
creates a sparse matrix made of $n \times m$ diagonal blocks.

\section{Eliminators}

The suggested method requires also the use of
 a certain type of matrices called \textit{eliminators} \cite{Podlubny-MFC-2000},
which are obtained from the $N\times N$ unit matrix $E$ by keeping only some
of its rows and omitting all other rows:
$S_1$ is obtained by omitting only the first row of $E$;
$S_2$ is obtained by omitting only the second row;
$S_{1,2}$ is obtained by omitting only the first and the second row of $E$;
and, in general,
$S_{r_1, r_2, \ldots, r_k}$ is obtained by omitting the rows
with the numbers $r_1, r_2, \ldots, r_k$.
In case of infinite matrices, similar matrices appeared in \cite{Cooke-IMS-book}.

If $A$ is a square $N\times N$ matrix, then
the product $S_{r_1, r_2, \ldots, r_k}A$
contains only rows of $A$ with the numbers
different from $r_1, r_2, \ldots, r_k$.
Similarly, the product $AS_{r_1, r_2, \ldots, r_k}^{T}$
contains only columns of $A$ with the numbers
different from $r_1, r_2, \ldots, r_k$.

The following simple example illustrates the main property of eliminators:
\begin{displaymath}
    A =
    \left[
    \begin{array}{lll}
    a_{11} & a_{12} & a_{13} \\
    a_{21} & a_{22} & a_{23} \\
    a_{31} & a_{32} & a_{33} \\
    \end{array}
    \right];
    \qquad
    S_{1} =
    \left[
    \begin{array}{lll}
    0 & 1 & 0 \\
    0 & 0 & 1 \\
    \end{array}
    \right];
    \qquad
    S_{1} \, A =
    \left[
    \begin{array}{lll}
    a_{21} & a_{22} & a_{23} \\
    a_{31} & a_{32} & a_{33} \\
    \end{array}
    \right];
\end{displaymath}
\begin{displaymath}
    A \, S_{1}^{T} =
    \left[
    \begin{array}{ll}
    a_{12} &  a_{13} \\
    a_{22} &  a_{23} \\
    a_{32} &  a_{33} \\
    \end{array}
    \right];
    \qquad
    S_{1} \, A \, S_{1}^{T} =
    \left[
    \begin{array}{ll}
    a_{22} &  a_{23} \\
    a_{32} &  a_{33} \\
    \end{array}
    \right].
\end{displaymath}

\section{Shifters}

For some types of approximation of differential operators 
(for example, one of the approximations of the symmetric Riesz derivative
below in this article) and especially for 
numerical solution of differential equations of arbitrary order (integer or fractional)
with delays,
it is convenient to introduce another special kind of matrices -- \emph{shifters} --,
which will represent discrete shifts, like, for example, delays.

Shifters (although without using this term) were used in \cite{Podlubny-MFC-2000}
for a simple generation of triangular strip matrices. 
There are shifters of two kinds:
$(N+1) \times (N+1)$ matrices
$E_{N,p}^{+}$, $p=1, \ldots N$,  with ones on $p$-th  diagonal
above the main diagonal and zeroes elsewhere,
and matrices $E_{N,p}^{-}$, $p=1, \ldots N$, with ones on $p$-th diagonal
below the main diagonal and zeroes elsewhere.
We also denote $E_{N,0}^{\pm} \equiv E_{N}$ the unit matrix.

The shift of all the coefficients in the triangular strip matrix $U_N$
in the south-west (bottom-left) direction can be easily written 
if we start with $U_{N+1}$ and then use shifters and eliminators:

\begin{equation}\label{eq:shift-sw}
_{-1}U_{N}
= 
S_1 \, E_{N+1,1}^{-} \, U_{N+1} \, E_{N+1,1}^{-} \, S_{N+1}^T 
\end{equation}

Similarly, the shift of all the coefficients in the triangular strip matrix $U_N$
in the north-east (top-right) direction can be easily obtained as:

\begin{equation}\label{eq:shift-ne}
_{+1}U_{N}
= 
S_{N+1} \, E_{N+1,1}^{+} \, U_{N+1} \, E_{N+1,1}^{+} \, S_{1}^T 
\end{equation}

\section{Discretization of ordinary fractional derivatives}

It follows from \cite{Podlubny-MFC-2000}, that the left-sided Riemann-Liouville or Caputo fractional derivative $v^{(\alpha)} (t) = \, _{0}D_{t}^{\alpha}v(t)$ can be approximated in all nodes of the equidistant discretization net 
$t = j\tau$  ($j = 0, 1, \ldots, n$) simultaneously with the help of the upper triangular strip 
matrix $B_n^{(\alpha)}$ as
\footnote{
In this article due to the use of the descending numbering of discretization nodes the roles 
of the matrices $B_n^{(\alpha)}$ (originally for backward fractional differences) 
and $F_n^{(\alpha)}$ (originally for forward fractional differences) are swapped
in comparison with  \cite{Podlubny-MFC-2000}, where these matrices were introduced
for the first time. However,  we would like to preserve the notation $B_n^{(\alpha)}$
for the case of the backward fractional differences approximation and $F_n^{(\alpha)}$
for the case of the forward fractional differences approximation. 
}
:

\begin{equation}\label{eq:left-sided-approximation}
\Bigl[
v_n^{(\alpha)} \;\;  
v_{n-1}^{(\alpha)} \;\; 
 \ldots \; 
 v_{1}^{(\alpha)} \;\; 
 v_{0}^{(\alpha)}
\Bigr]^T
=
 B_{n}^{(\alpha)} \;
\Bigl[
v_n \;\;  v_{n-1} \;\;  \ldots \; v_{1} \;\; v_{0}
\Bigr]^T
\end{equation}

\noindent
where
\begin{equation}
    B_{n}^{(\alpha)}
    =
     \frac{1}{\tau^\alpha}
    \left[
    \!\!
    \begin{array}{cccccc}
    \omega^{(\alpha)}_0     & \omega^{(\alpha)}_1       & \ddots                & \ddots                & \omega^{(\alpha)}_{n-1} & \omega^{(\alpha)}_n \\
    0                       & \omega^{(\alpha)}_0       & \omega^{(\alpha)}_1   & \ddots                & \ddots                & \omega^{(\alpha)}_{n-1} \\
    0                       & 0                         & \omega^{(\alpha)}_0   & \omega^{(\alpha)}_1   & \ddots                & \ddots                \\
    \cdots                  & \cdots                    & \cdots                & \ddots                & \ddots                & \ddots                \\
    0                       & \cdots                    & 0                     & 0                     & \omega^{(\alpha)}_0   & \omega^{(\alpha)}_1   \\
    0                       & 0                         & \cdots                & 0                     & 0                     & \omega^{(\alpha)}_0   \\ 
    \end{array}
    \!\!\!
    \right]
\end{equation}

\begin{equation}\label{eq:Omega-Coefficients-definition}
    \omega_j^{(\alpha)} = (-1)^j { \alpha \choose j},
    \qquad
    j = 0, 1, \ldots, n.
\end{equation}

Similarly, the right-sided Riemann-Liouville or Caputo fractional derivative $v^{(\alpha)} (t) = \, _{t}D_{b}^{\alpha}v(t)$ can be approximated in all nodes of the equidistant discretization net 
$t = j \tau$  ($j = 0, 1, \ldots, n$) simultaneously with the help of the lower triangular strip 
matrix $F_n^{(\alpha)}$:

\begin{equation}\label{eq:right-sided-approximation}
\Bigl[
v_n^{(\alpha)} \;\;  
v_{n-1}^{(\alpha)} \;\; 
 \ldots \; 
 v_{1}^{(\alpha)} \;\; 
 v_{0}^{(\alpha)}
\Bigr]^T
=
 F_{n}^{(\alpha)} \; 
\Bigl[
v_n \;\;  v_{n-1} \;\;  \ldots \; v_{1} \;\; v_{0}
\Bigr]^T
\end{equation}

\begin{equation}
F_n^{(\alpha)}
  =
    \frac{1}{\tau^\alpha}
    \left[
    \!\!
    \begin{array}{cccccc}
    \omega^{(\alpha)}_0   & 0  & 0  & 0  & \cdots & 0 \\
    \omega^{(\alpha)}_1   & \omega^{(\alpha)}_0   & 0  & 0  & \cdots & 0 \\
    \omega^{(\alpha)}_2   & \omega^{(\alpha)}_1  & \omega^{(\alpha)}_0 & 0  & \cdots & 0 \\
    \ddots & \ddots & \ddots & \ddots & \cdots & \cdots \\
    \omega^{(\alpha)}_{n-1}    & \ddots  & \omega^{(\alpha)}_2   & \omega^{(\alpha)}_1 & \omega^{(\alpha)}_0  & 0 \\
    \omega^{(\alpha)}_n  & \omega^{(\alpha)}_{n-1}   & \ddots  & \omega^{(\alpha)}_2  & \omega^{(\alpha)}_1  & \omega^{(\alpha)}_0  \\
    \end{array}
    \!\!\!
    \right]
\end{equation}

The symmetric Riesz derivative of order $\beta$ can be approximated based on its definition
(\ref{eq:Riesz-derivative-definition}) as a combination of the approximations 
(\ref{eq:left-sided-approximation}) and (\ref{eq:right-sided-approximation})
for the left-  and right-sided Riemann-Liouville derivatives, or using the centred fractional differences approximation of the symmetric Riesz derivative suggested recently by Ortigueira \cite{Ortigueira-2006,Ortigueira-2008}. The general formula is the same:

\begin{equation}\label{eq:symmetric-approximation}
\Bigl[
v_m^{(\beta)} \;\;  
v_{m-1}^{(\beta)} \;\; 
 \ldots \; 
 v_{1}^{(\beta)} \;\; 
 v_{0}^{(\beta)}
\Bigr]^T
=
 R_{m}^{(\beta)} \; 
\Bigl[
v_m \;\;  v_{m-1} \;\;  \ldots \; v_{1} \;\; v_{0}
\Bigr]^T
\end{equation}

In the first case, the approximation for the left-sided Caputo derivative is taken 
one step ahead, and the approximation for the right-sided Caputo derivatve
is taken one step back.  This leads to the matrix 

\begin{equation}\label{eq:ran-half-sum}
 R_{m}^{(\beta)} = 
 \frac{h^{-\alpha}}{2}
 \Bigl[ \; 
 	_{-1}U_{m} + \, _{+1}U_{m}
 \Bigr]
\end{equation}

In the second case (Ortigueira's definition \cite{Ortigueira-2006}) we have the following symmetric matix:

\begin{equation}\label{eq:ran-ort}
 R_{m}^{(\beta)} = 
 h^{-\beta} 
 \left[ \; 
    \!\!
    \begin{array}{cccccc}
    \omega^{(\beta)}_0   & \omega^{(\beta)}_1  & \omega^{(\beta)}_2  & \omega^{(\beta)}_3  & \cdots & \omega^{(\beta)}_m \\
    \omega^{(\beta)}_1   & \omega^{(\beta)}_0   & \omega^{(\beta)}_1  & \omega^{(\beta)}_2  & \cdots & \omega^{(\beta)}_{m-1} \\
    \omega^{(\beta)}_2   & \omega^{(\beta)}_1  & \omega^{(\beta)}_0 & \omega^{(\beta)}_1  & \cdots & \omega^{(\beta)}_{m-2} \\
    \ddots & \ddots & \ddots & \ddots & \cdots & \cdots \\
    \omega^{(\beta)}_{m-1}    & \ddots  & \omega^{(\beta)}_2   & \omega^{(\beta)}_1 & \omega^{(\beta)}_0  & \omega^{(\beta)}_1 \\
    \omega^{(\beta)}_m  & \omega^{(\beta)}_{m-1}   & \ddots  & \omega^{(\beta)}_2  & \omega^{(\beta)}_1  & \omega^{(\beta)}_0  \\
    \end{array}
    \!\!\!
 \right]
\end{equation}

\medskip

\begin{equation}
\omega_k^{(\beta)}
= 
\frac{(-1)^k \, \Gamma(\beta + 1) \,  \cos (\beta \pi /2)}
	{\Gamma(\beta/2 - k + 1) \, \Gamma(\beta/2 + k + 1)}, 
	\qquad 
	k = 0, 1, \ldots, m
\end{equation}

Both these approximations of symmetric Riesz derivatives give practically the same 
numerical results and in case of numerical solution of partial fractional differential equations
lead to a well-posed matrix of the resulting algebraic system.

\section{Discretization of partial derivatives in time and space}

\begin{figure}% 
	\centering 
	\begin{minipage}[t]{4cm}% 
		\centering
		\includegraphics[width=3.4cm]{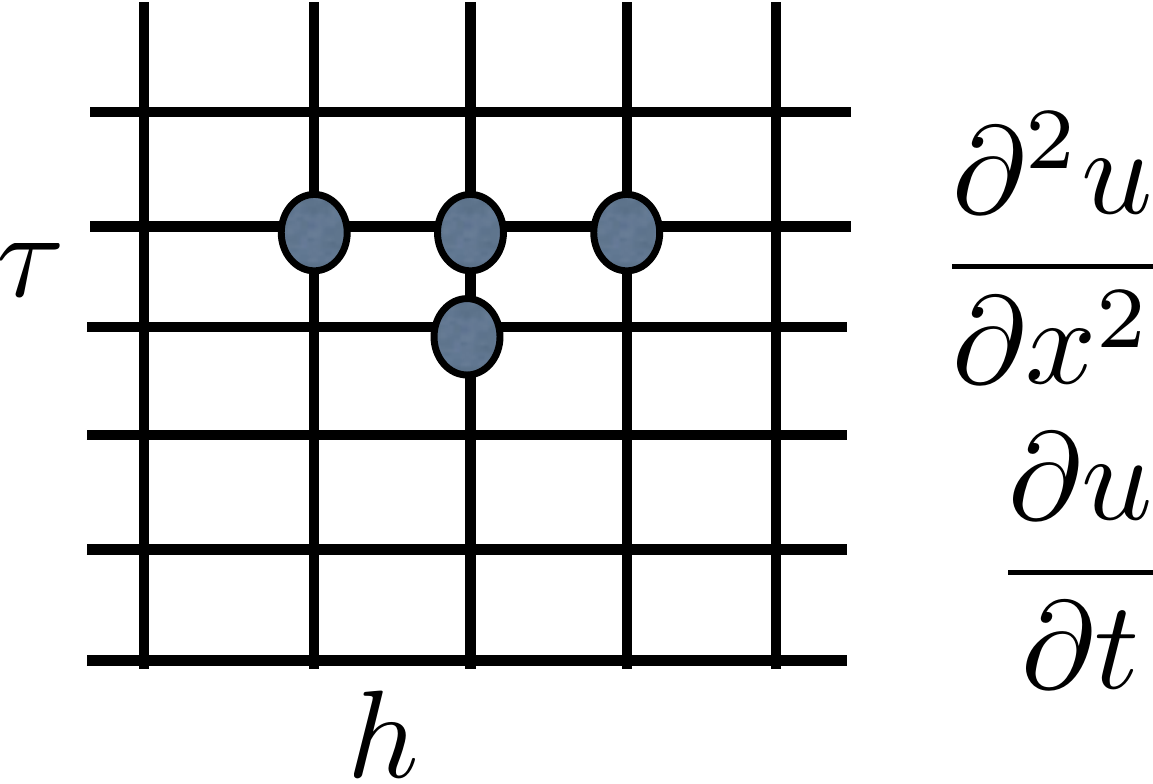} 
		\caption{A stencil for integer-order derivatives.}% 
		\label{fig:stencil-1}% 
	\end{minipage}% 
	\hfill
	\begin{minipage}[t]{4cm}% 
		\centering
		\includegraphics[width=3.7cm]{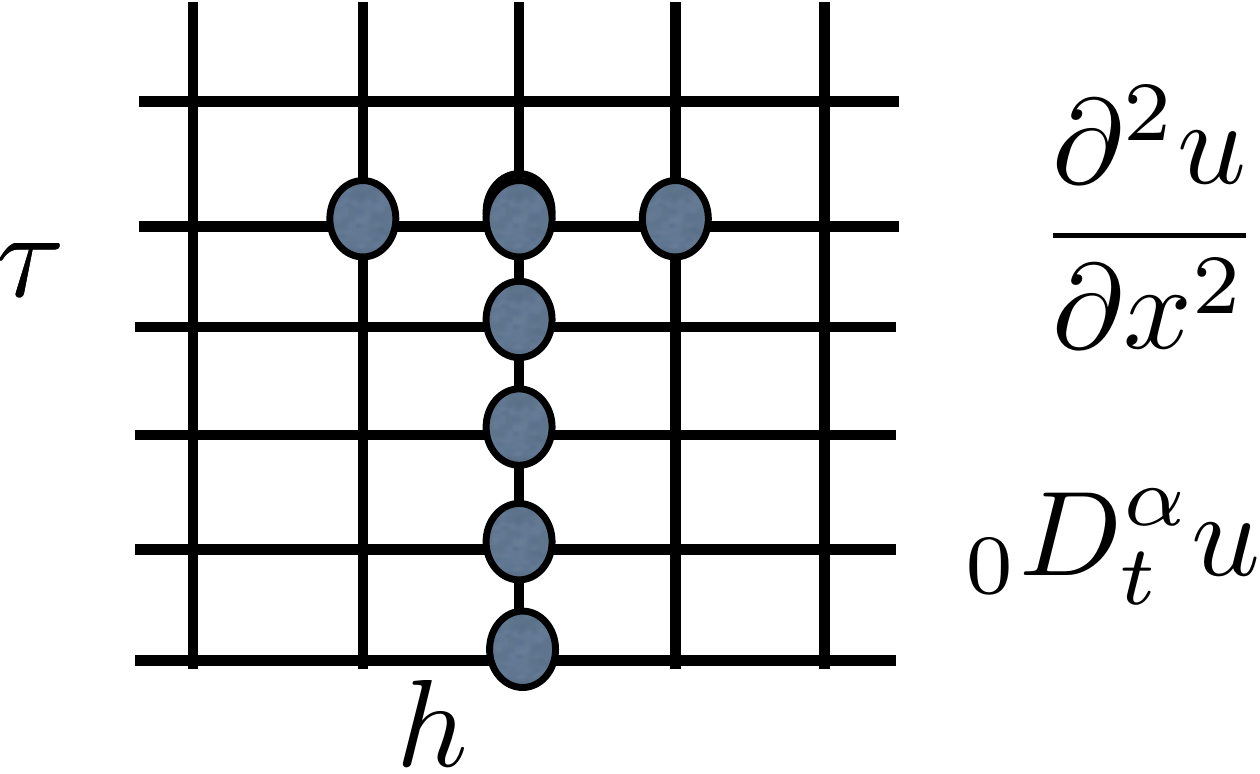} 
		\caption{A stencil in case of fractional time derivative.}% 
		\label{fig:stencil-2}% 
	\end{minipage}% 
	\hfill
	\begin{minipage}[t]{4cm}% 
		\includegraphics[width=3.7cm]{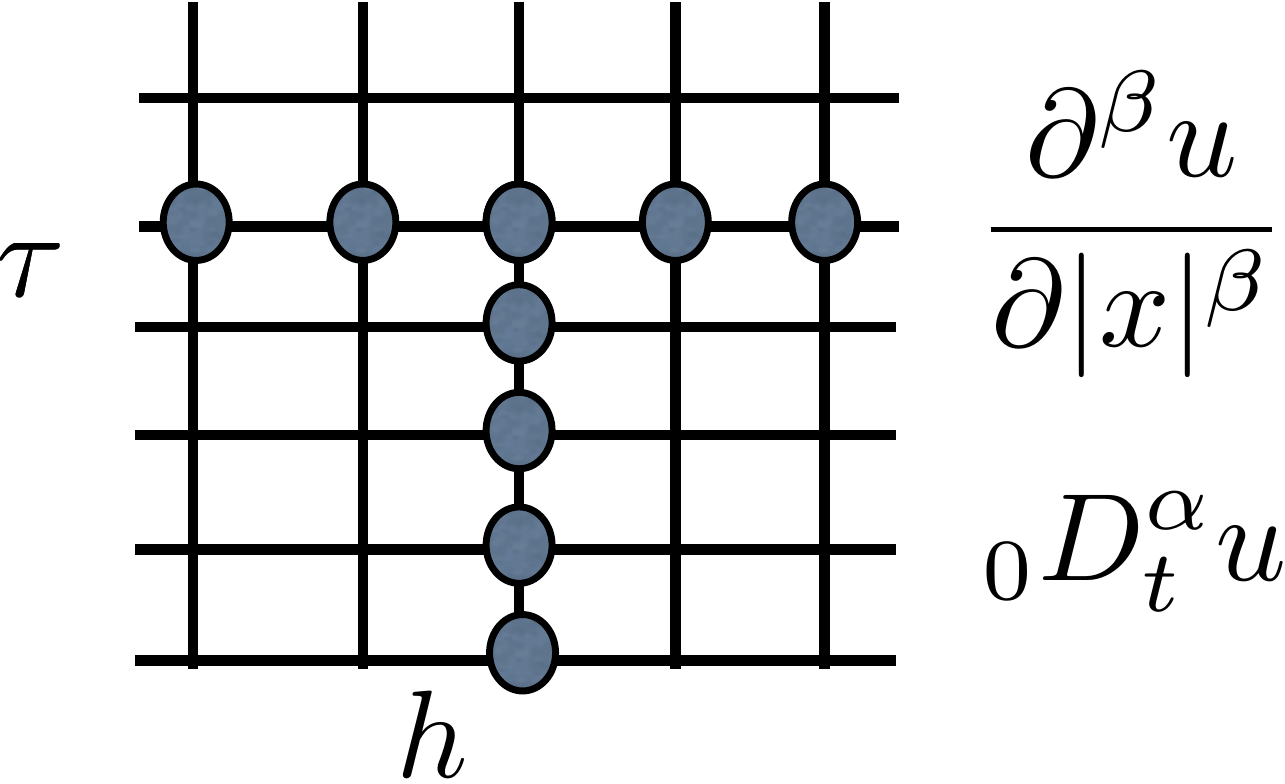} 
		\caption{A stencil in case of fractional time and spatial derivatives.}% 
		\label{fig:stencil-3}% 
	\end{minipage}% 
\end{figure}

The simplest implicit discretization scheme for the classical diffusion equation is shown in Fig.~\ref{fig:stencil-1},
where the two nodes in time direction are used for approximating the time derivative, 
and the three points in spatial direction are used for the symmetric approximation 
of the  the spatial derivative.
The stencil in  Fig.~\ref{fig:stencil-1} involves therefore only two time layers.
If we consider fractional-order time derivative, then we have to involve all 
time levels starting from the very beginning. This is shown in Fig.~\ref{fig:stencil-2}
for the case of five time layers.

Similarly, if in addition to fractional-order time derivative we also consider 
symmetric fractional-order spatial derivatives, 
then we have to use all nodes at the considered time layer
from the leftmost to the rightmost spatial discretization node.
This most general situation is shown in  Fig.~\ref{fig:stencil-3}.

Let us consider the nodes $(ih, j\tau)$, $j = 0, 1, 2, \ldots, n$, 
corresponding to all time layers at $i$-th spatial discretization node. 
It has been shown in \cite{Podlubny-MFC-2000} that all values 
of  $\alpha$-th order time derivative of $u(x,t)$ at these nodes are approximated
using the discrete analogue of differentiation of arbitrary order:

\begin{equation}
\Bigl[u_{i,n}^{(\alpha)} \; u_{i,n-1}^{(\alpha)} \; \ldots \; u_{i,2}^{(\alpha)} \; u_{i,1}^{(\alpha)} \; u_{i,0}^{(\alpha)}\Bigr]
=
B_n^{(\alpha)} \, 
\Bigl[u_{i,n} \; u_{i,n-1}  \; \ldots \; u_{i,2} \; u_{i,1} \; u_{i,0}\Bigr]^T
\end{equation}

In order to obtain a simultaneous approximation of  $\alpha$-th order time derivative of $u(x,t)$
in all nodes shown in Fig.~\ref{fig:nodes-net}, we need to arrange all function values $u_{ij}$ at the discretization nodes to the form of a column vector:

\begin{eqnarray}
u_{nm} = 
\Bigl[
u_{m,n} \; u_{m-1, n} \;  \ldots \; u_{1,n} \; u_{0, n} \;  \hspace*{8em} && \nonumber \\
u_{m,n-1} \; u_{m-1, n-1} \;  \ldots \; u_{1,n-1} \; u_{0, n-1} \; \hspace*{2.5em}  && \nonumber\\
\ldots  \ldots \ldots  \hspace*{11em} && \nonumber\\
u_{m,1} \; u_{m-1, 1} \;  \ldots \; u_{1,1} \; u_{0, 1} \; \hspace*{3em}  && \nonumber\\
u_{m,0} \; u_{m-1, 0} \;  \ldots \; u_{1,0} \; u_{0, 0} \;  
\Bigr]^T
\end{eqnarray}

In visual terms of Fig.~\ref{fig:nodes-net}, we first take the nodes of $n$-th time layer, 
then the nodes of $(n-1)$-th time layer, and so forth, and put them in this order in a vertical column stack. 

The matrix that transforms the vector $U_{nm}$ to the vector $U_{t}^{(\alpha)}$ of the partial fractional
derivative of order $\alpha$ with respect to time variable can be obtained 
as a Kronecker product of the matrix $B_{n}^{(\alpha)}$,  which  corresponds to 
the fractional ordinary derivative of order $\alpha$ (recall that $n$ is the number of time steps), 
and the unit matrix $E_{m}$ (recall that $m$ is the number of spatial discretization nodes):

\begin{equation}
T_{mn}^{(\alpha)} = B_{n}^{(\alpha)} \otimes E_{m}
\end{equation}

\begin{figure}% 
	\centering 
	\begin{minipage}[t]{6cm}% 
		\centering
		\includegraphics[width=5.5cm]{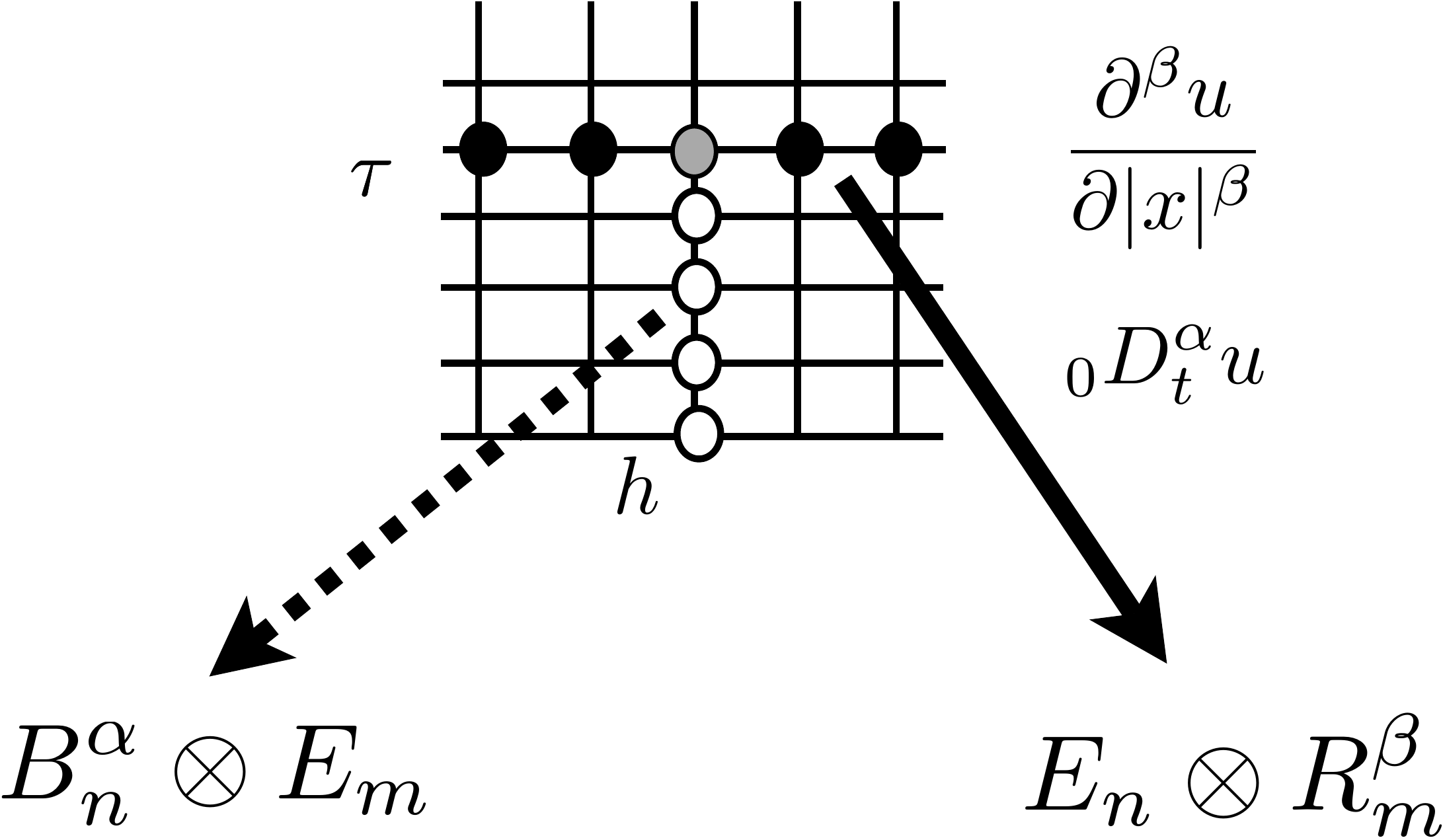} 
		\caption{Discretization of partial derivatives.}% 
		\label{fig:discretization-partial-derivatives}% 
	\end{minipage}% 
	\hfill
	\begin{minipage}[t]{6cm}% 
		\centering
		\includegraphics[width=6.5cm]{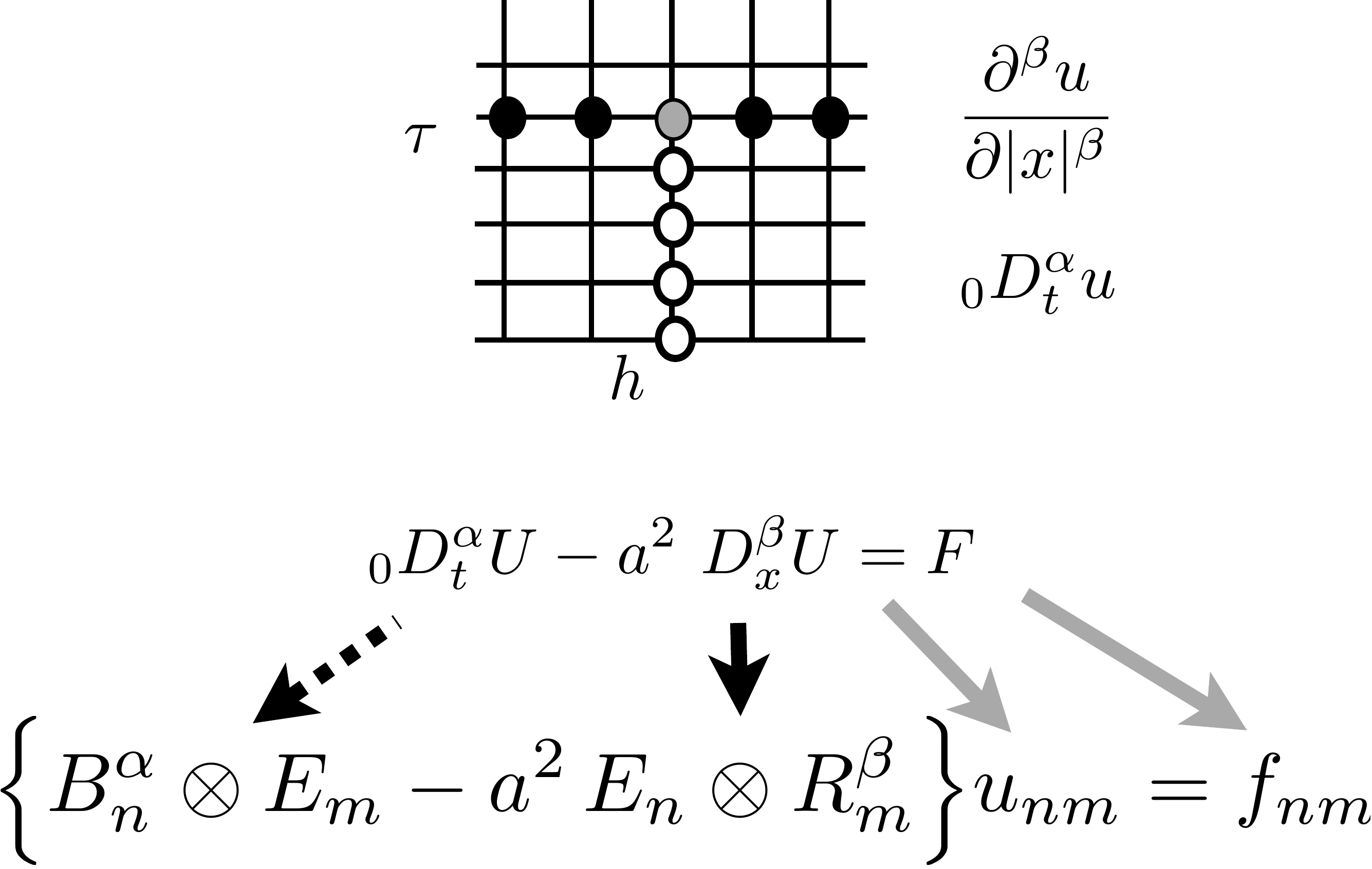} 
		\caption{Discretization of partial derivatives and of the equation}% 
		\label{fig:discretization-equation}% 
	\end{minipage}% 
\end{figure}

This is illustrated in Fig.~\ref{fig:discretization-partial-derivatives}, where the
nodes denoted as white and gray are used to approximate the fractional-oder
time derivative at the node shown in gray.

Similarly, the matrix that transforms the vector $U$ to the vector $U_{x}^{(\beta)}$ of the partial fractional derivative of order $\beta$ with respect to spatial variable can be obtained 
as a Kronecker product of 
the unit matrix $E_{n}$ (recall that $n$ is the number of spatial discretization nodes),
and the matrix $R_{m}^{(\beta)}$,  which  corresponds to a symmetric Riesz
ordinary derivative of order $\beta$
\cite{Ortigueira-2006,Ortigueira-2008}
 (recall that $m$ is the number of time steps):

\begin{equation}
S_{mn}^{(\alpha)} = E_{n}  \otimes  R_{n}^{(\beta)}
\end{equation}

This is also illustrated in Fig.~\ref{fig:discretization-partial-derivatives}, where the
nodes denoted as black and gray (corresponding to all discretization nodes from the leftmost to the rightmost one) are used to approximate the symmetric fractional-order
Riesz derivative at the same node shown in gray.

\begin{figure}
\begin{center}
\includegraphics[width=8cm]{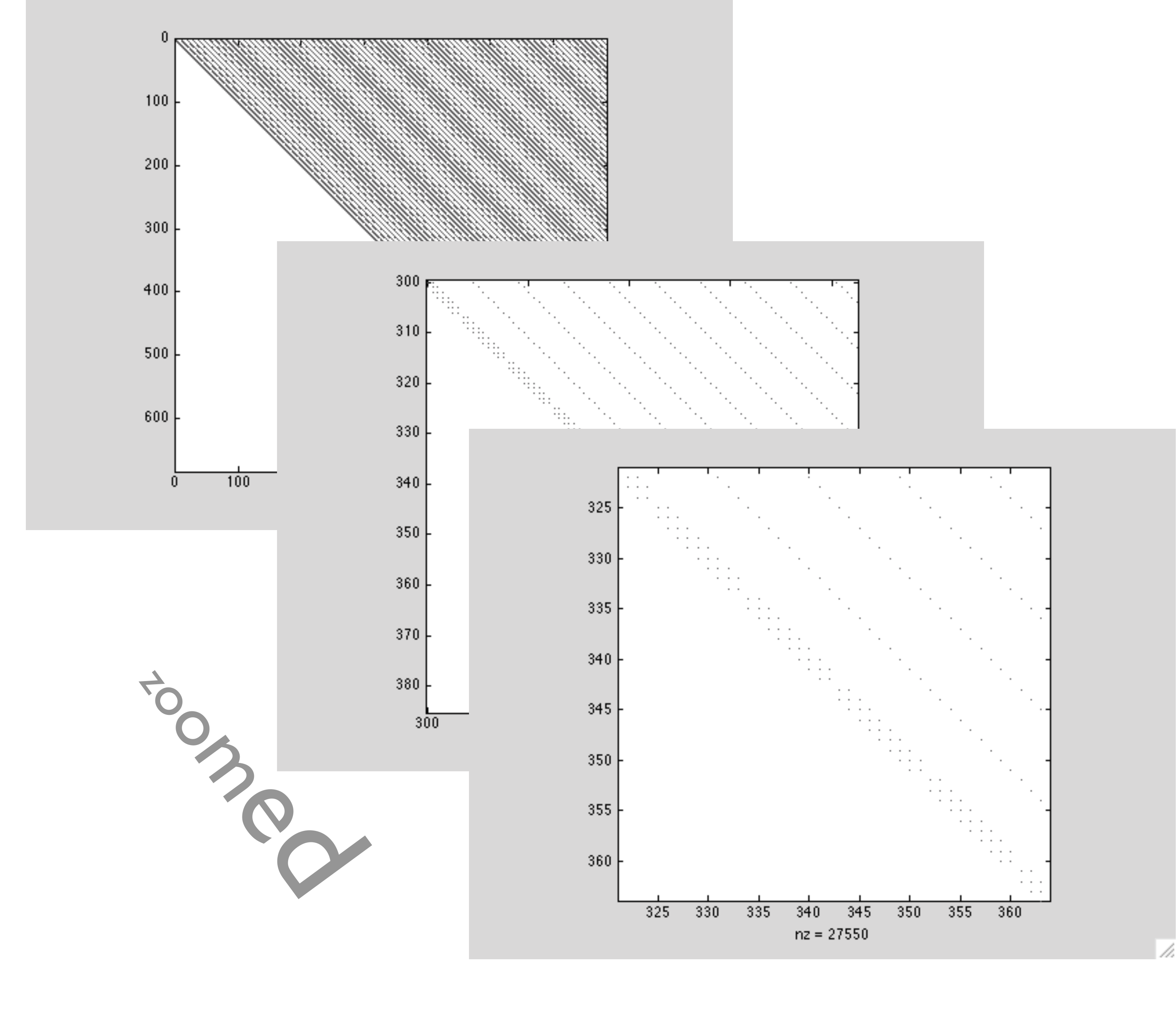}
\caption{The structure of the matrix of the resulting algebraic system.}
\label{fig:matrix-structure-zoomed}
\end{center}
\end{figure}

Having these approximations for partial fractional derivatives with respect to both variables, 
we can immediately discretize the general form of the fractional diffusion equation 
by simply replacing  the derivatives with their discrete analogs (Fig.~\ref{fig:discretization-equation}). Namely, the equation 

\begin{equation}
_{0}^{C}D_{t}^\alpha u -
\chi 
\frac{\partial^\beta u}{\partial |x|^\beta}
= f(x,t)
\end{equation}

\noindent
is discretized as 

\begin{equation}
\Bigl\{
B_{n}^{(\alpha)} \otimes E_m - \chi \,
E_n \otimes R_{m}^{(\beta)}
\Bigr\} u_{nm} = f_{nm},
\end{equation}

\noindent
and the matrix of this system has the structure shown in Fig.~\ref{fig:matrix-structure-zoomed}.

\section{Initial and boundary conditions}

Initial and boundary conditions must be equal to zero.
If it is not so, then an auxiliary unknown function must be introduced,
which satisfies the zero initial and boundary conditions. 
In this way, the non-zero initial and boundary conditions moves
to the right-hand side of the equation for the new unknown function.

\section{Implementation in MATLAB}

We provide a set of MATLAB routines for implementing the suggested method \cite{MFC2demo}.
The function BCRECUR returns the values  of the coefficients that appear 
in the fractional difference approximations of fractional derivatives. 
The function BAN returns the matrix for the backward difference approximation
of the left-sided ordinary fractional derivative, the function FAN returns the matrix 
for approximating the right-sided ordinary fractional derivative,
and the functions RANSYM and RANORT return the matrices for 
approximating the symmetric Riesz using the formulas
(\ref{eq:ran-half-sum}) and (\ref{eq:ran-ort}), respectively.
The function ELIMINATOR returns the eliminator matrix,
and the function SHIFT implements the operations (\ref{eq:shift-sw}) and (\ref{eq:shift-ne}).

The use of these routines is illustrated by the demo functions
FRAC\-DIFF\-DEMOU, which implements Examples 1 and 2 below,
 FRAC\-DIFF\-DEMOY, which implements Examples 3 and 4,
 and FRAC\-DIFF\-DEMOY\-DELAY, which implements Example 5.

\section{Examples}

In  this section we introduce several examples illustrating the use of the suggested method. 

First, we demonstrate that for the classical integer-order diffusion equation
our method gives proper results, which are in agreement with the analytical 
and numerical solutions provided in \cite{milne}. 

Second, we obtain the numerical solution of a time-fractional diffusion equation.
This solution is in perfect agreement with the numerical solution obtained
in the very recent work \cite{Scherer:2008kk} by a different approach. 

Then we consider fractional diffusion equation with spatial fractional derivative.
The fractional derivative with respect to the spatial variable is  considered 
as a Riesz fractional derivative.

After that, we show the results of numerical solution of a general fractional 
diffusion equation, where time and spatial derivatives are both of fractional order --
the time fractional derivative is a left-sided Riemann--Liouville derivative, 
and the spatial fractional derivative is a Riesz fractional derivatie. 

Finally, we demonstrate that consideration of partial differential equations with 
fractional derivatives and delays is equally simple in the framework of the 
suggested general approach.

In all examples, the spatial interval is finite.

\subsection{Example 1: Classical diffusion equation}

Let us start with the classical problem \cite{milne}:

\begin{equation} \label{eq:example-1-equation}
\frac{\partial u}{\partial t} =
\frac{\partial^2 u}{\partial x^2}
\end{equation}

\begin{equation} \label{eq:example-1-bc}
u(0,t) = 0, \quad u(1,t) = 0
\end{equation}

\begin{equation} \label{eq:example-1-ic}
u(x,0) = 4x (1-x)
\end{equation}

To reduce this problem to a problem with zero initial conditions (the boundary conditions are already zero), let us introduce an auxiliary function

\begin{equation} \label{eq:example-1-y-u}
y(x,t) = u(x,t) - u(x,0)
\end{equation}

It follows from (\ref{eq:example-1-y-u}) and (\ref{eq:example-1-equation})--(\ref{eq:example-1-ic})
that the function $y(x,t)$ must satisfy

\begin{equation} \label{eq:example-1-equation-y}
\frac{\partial y}{\partial t} -
\frac{\partial^2 y}{\partial x^2} = f(x,t),  \qquad (\mbox{with } f(x,t) \equiv  8)
\end{equation}

\begin{equation} \label{eq:example-1-ibc-y}
y(0,t) = 0, \quad y(1,t) = 0; 
\qquad  \quad
y(x,0) = 0.
\end{equation}

The problem (\ref{eq:example-1-equation-y})--(\ref{eq:example-1-ibc-y}) can be discretized using the described method (see Fig.~\ref{fig:discretization-equation}), which gives

\begin{equation}\label{eq:example-1-discretization}
\Bigl\{
B_{n}^{(1)} \otimes E_m - \,
E_n \otimes R_{m}^{(2)}
\Bigr\} y_{nm} = f_{nm}
\end{equation}

\noindent
where $m$ is the number of spatial discretization intervals and $n$ is the number of time steps.

To obtain the system for finding the unknown values of $y_{nm}$ for the inner nodes of the discretization net, we have to use the initial and boundary conditions. Since they all are zero, it is sufficient to delete the corresponding rows and columns in the system (\ref{eq:example-1-discretization}), which is easily done with the help of \textit{eliminators}.

\begin{figure}
\begin{center}
\includegraphics[width=6cm]{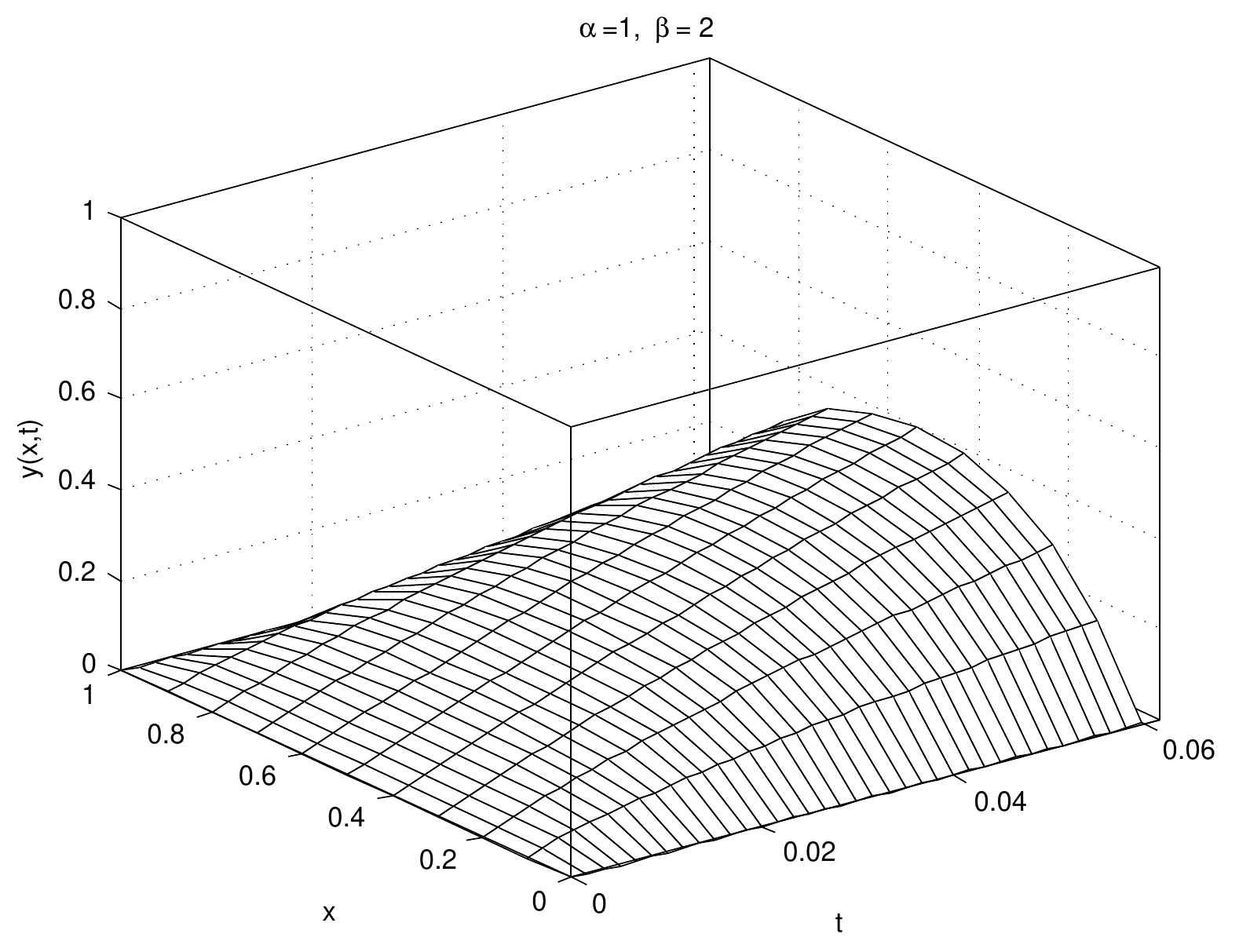} \hfill \includegraphics[width=6cm]{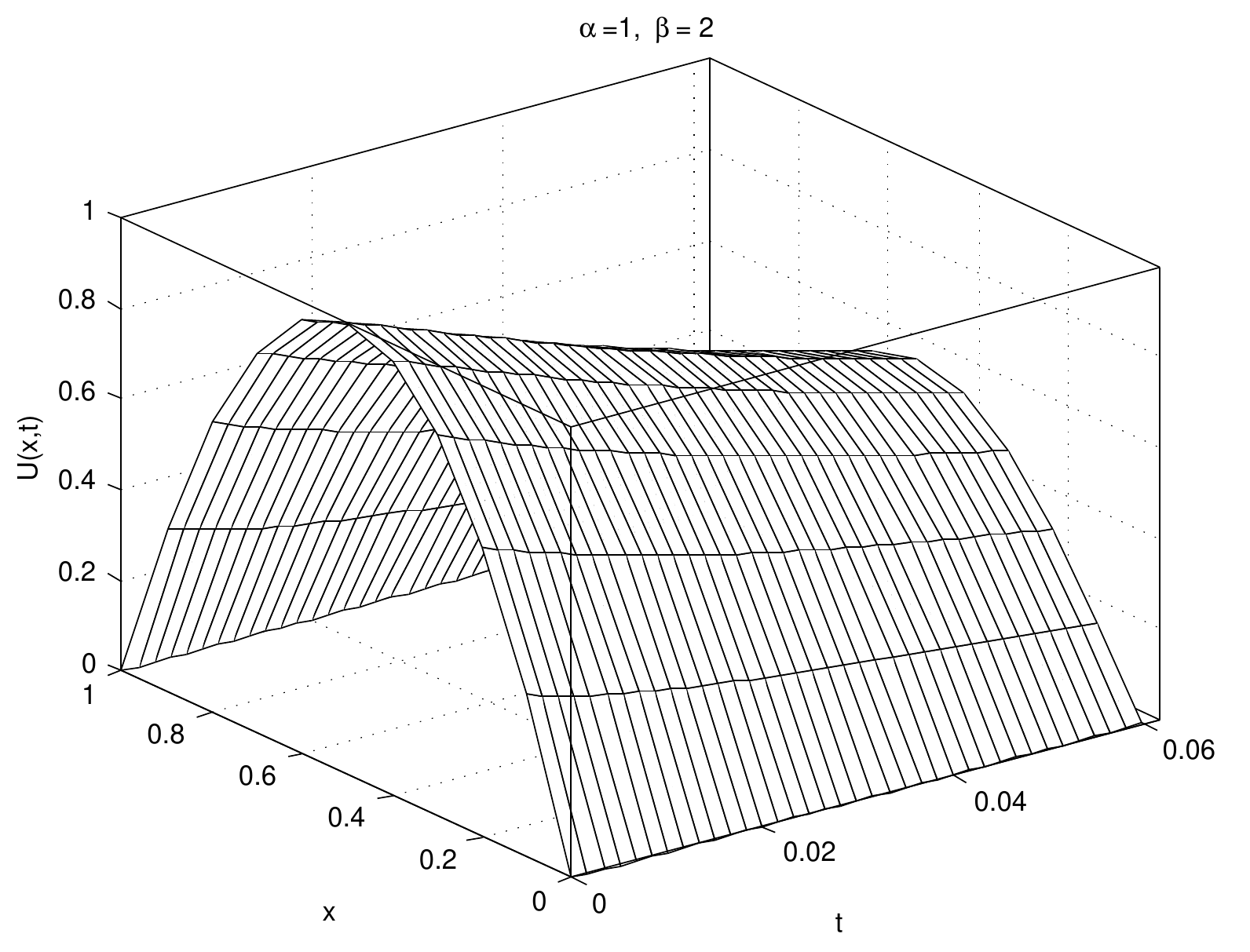}
\caption{Solutions $y(x,t)$ (left)  and $u(x,t)$ (right) of Example 1, with the same values of parameters as in \cite{milne}.}
\label{fig:example-1}
\end{center}
\end{figure}

The result of computation of $y(x,t)$ for the spatial step $h=0.1$ and the time step $\tau = h^2/6$
is shown in Fig.~\ref{fig:example-1} (on the left) for $n=37$ time steps. These values were chosen for the purpose of comparison with the results from \cite{milne}.
Using (\ref{eq:example-1-y-u}), we can compute $u(x,t)$, and the result is shown in Fig.~\ref{fig:example-1} (on the right). The values of $u(x,t)$ are in perfect agreement with the values given in \cite{milne} for the same values of $h$, $\tau$, and $n$.

\bigskip
\bigskip

\subsection{Example 2: Diffusion equation with time fractional derivative}

Now let us consider the equation with the Caputo fractional-order time derivative:

\begin{equation}\label{eq:example-2-equation}
 _{0}^{C}\!D_{t}^{\alpha}u =  \frac{\partial^2 u}{\partial x^2} 
\end{equation}

\begin{equation}\label{eq:example-2-bc}
u(0,t) = 0, \quad u(1,t) = 0
\end{equation}

\begin{equation} \label{eq:example-2-ic}
u(x,0) = 4x (1-x)
\end{equation}

Since the Caputo derivative of a constant is zero \cite{Caputo-1969,Podlubny-FDE-book}, 
for the auxiliary function $y(x,t)$ defined by equation (\ref{eq:example-1-y-u}) we obtain 
a problem with zero initial and boundary conditions similar to (\ref{eq:example-1-equation-y})--(\ref{eq:example-1-ibc-y}):

\begin{equation}
 _{0}^{C}\!D_{t}^{\alpha}y -  \frac{\partial^2 y}{\partial x^2} = f(x,t),  
  \qquad 
  (\mbox{with } f(x,t) \equiv  8) 
\end{equation}

\begin{equation}
y(0,t) = 0, \quad y(1,t) = 0; 
\qquad \quad 
y(x,0) = 0
\end{equation}
 
This problem can be discretized in the same manner as the previous one 
(refer to Fig.~\ref{fig:discretization-equation}), with the only difference 
that instead of the first-order time derivative we have now a derivative of order $\alpha$:

\begin{equation}\label{eq:example-2-discretization}
\Bigl\{
B_{n}^{(\alpha)} \otimes E_m - \,
E_n \otimes R_{m}^{(2)}
\Bigr\} y_{nm} = f_{nm}
\end{equation}

\noindent
where $m$ is the number of spatial discretization intervals and $n$ is the number of time steps.

As above, the use of the zero initial conditions means that the corresponding rows and columns in the system (\ref{eq:example-2-discretization}) are removed with the help of \textit{eliminators}.

\begin{figure}
\begin{center}
\includegraphics[width=6cm]{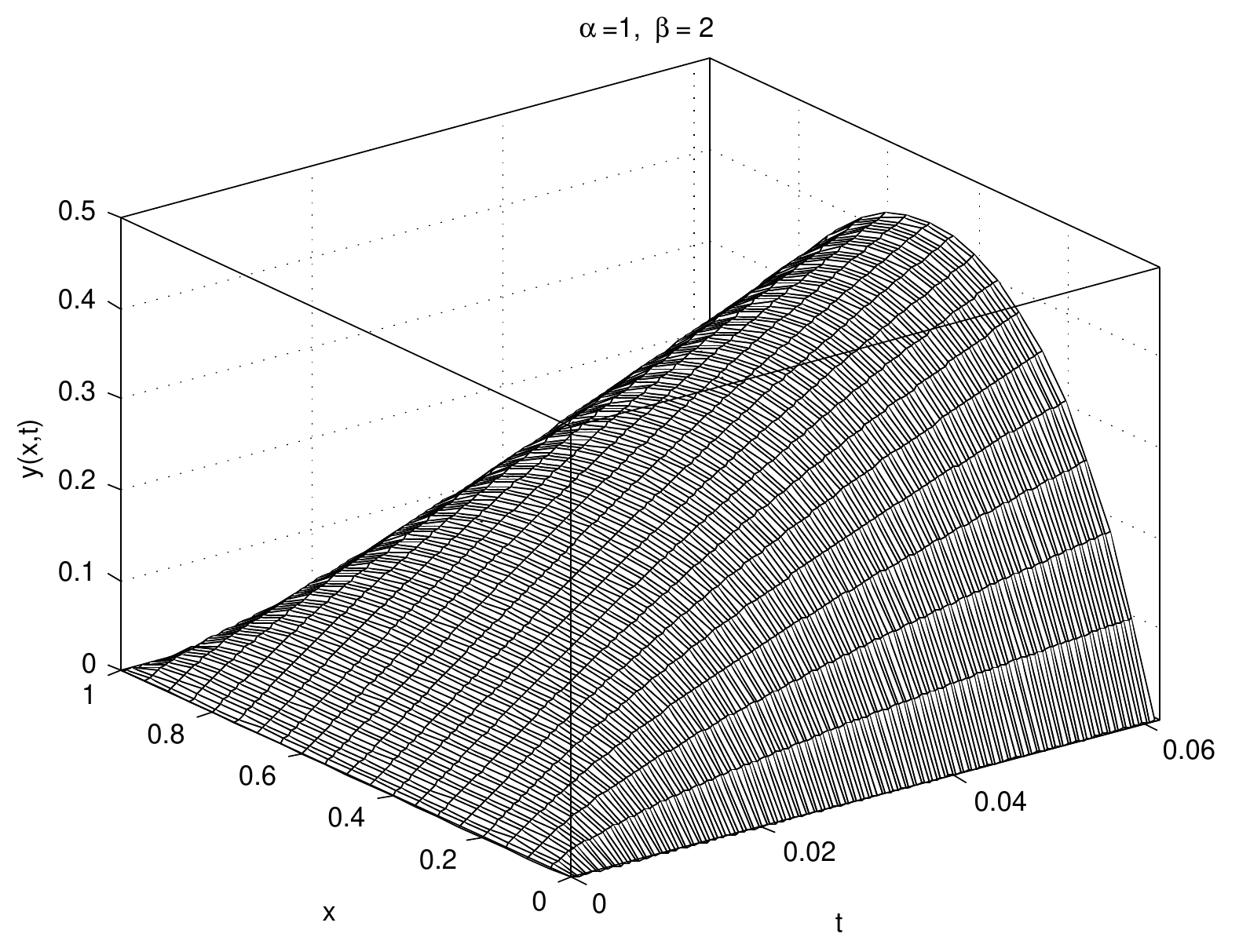} \hfill  \includegraphics[width=6cm]{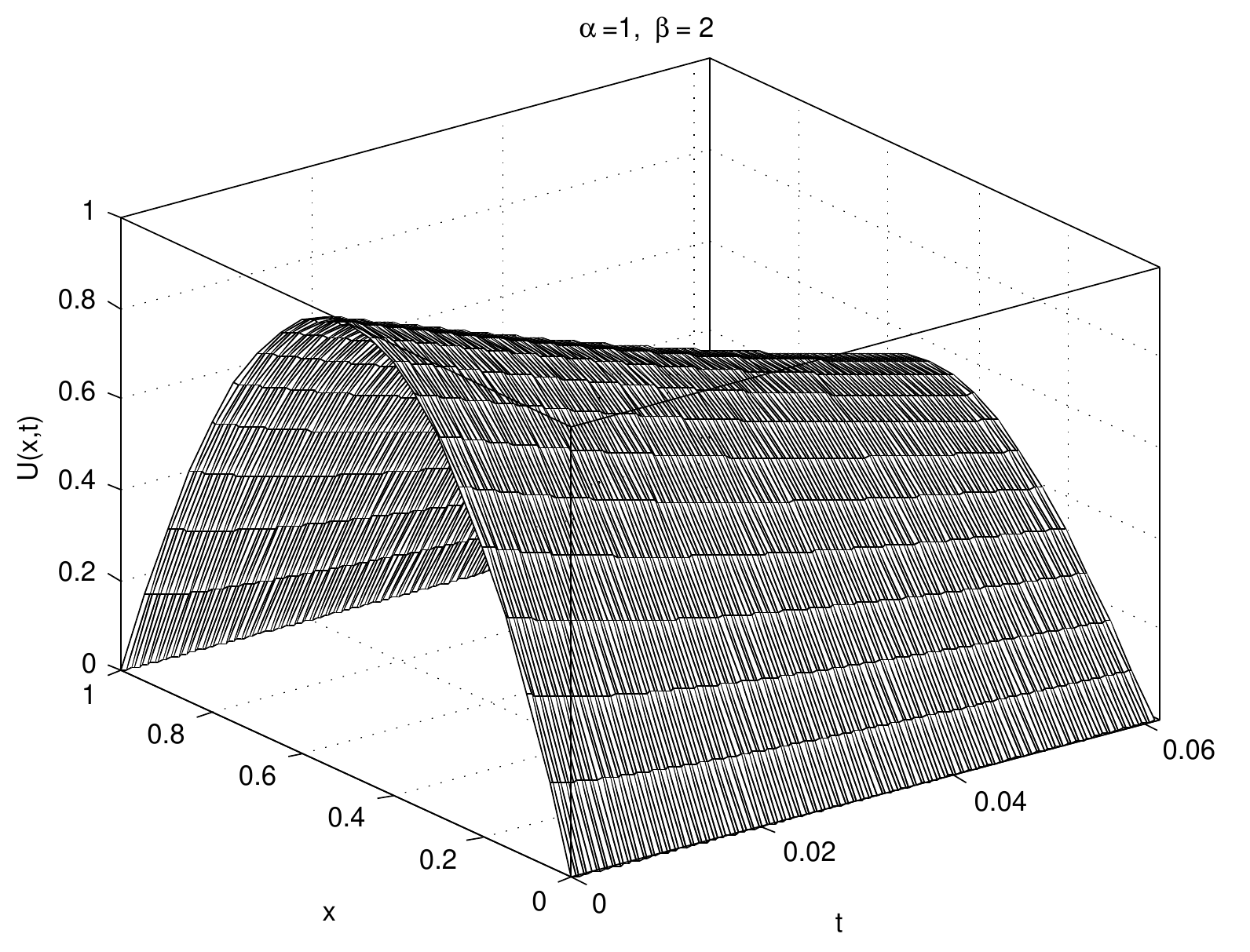}\\ 
\includegraphics[width=6cm]{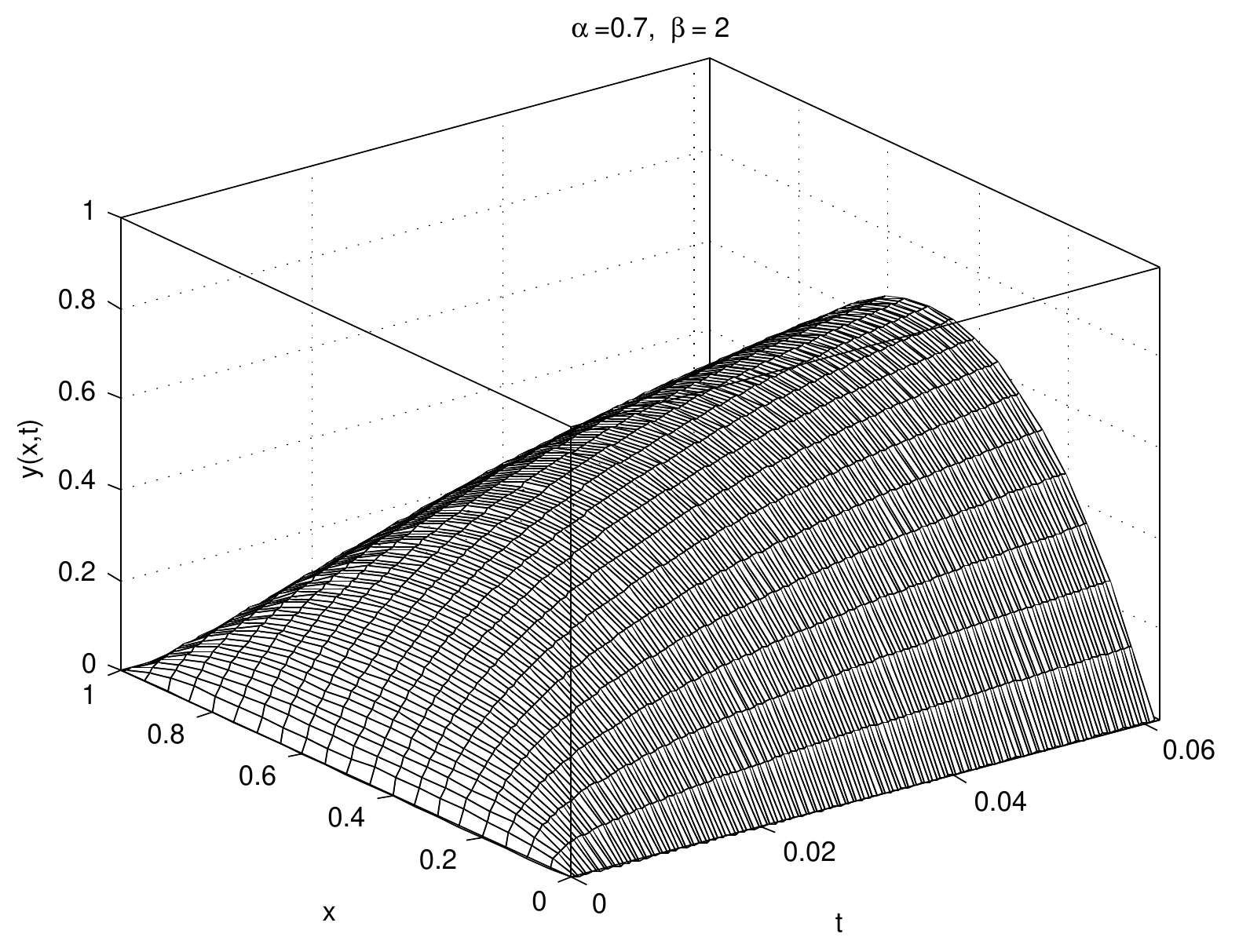} \hfill  \includegraphics[width=6cm]{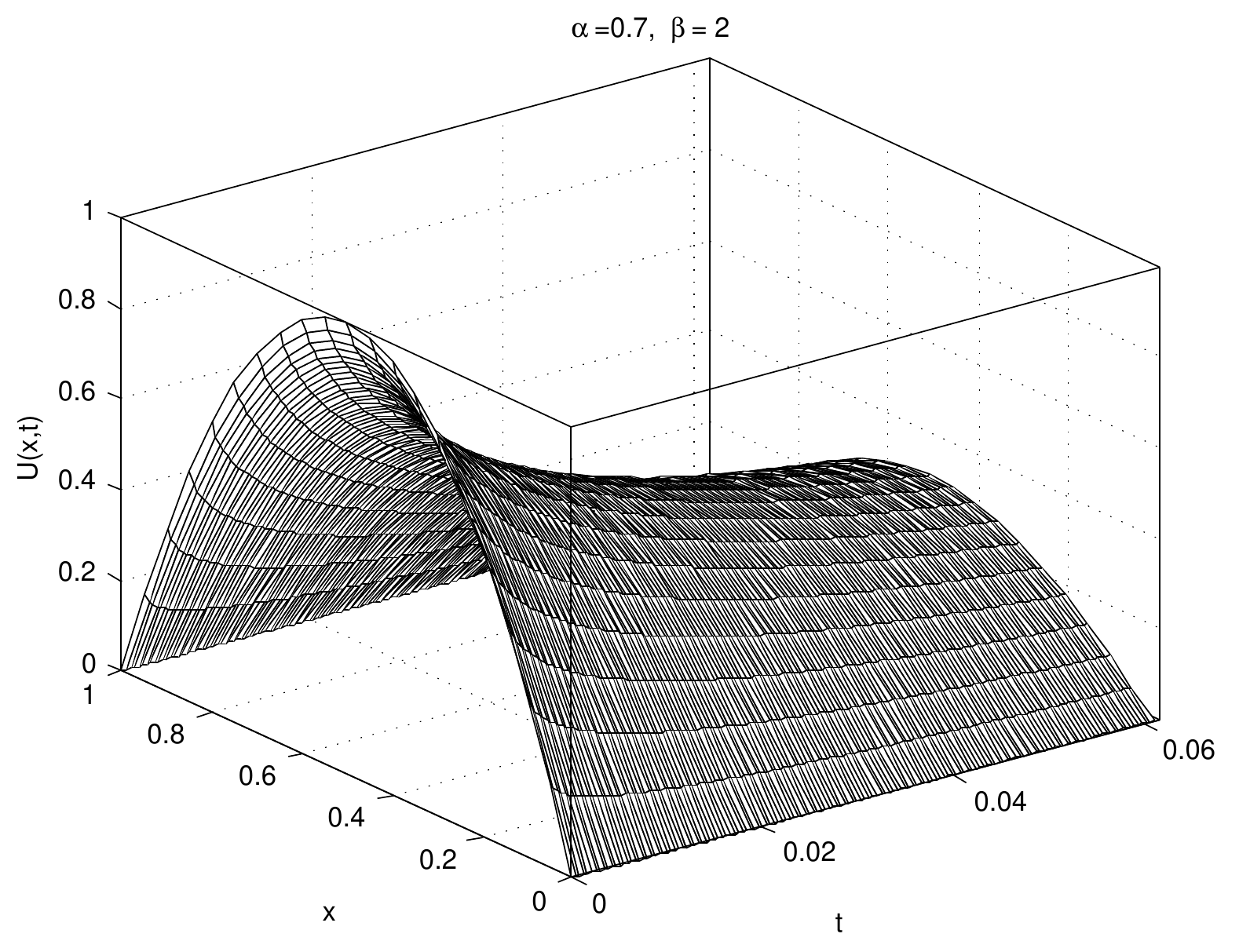}\\ 
\includegraphics[width=6cm]{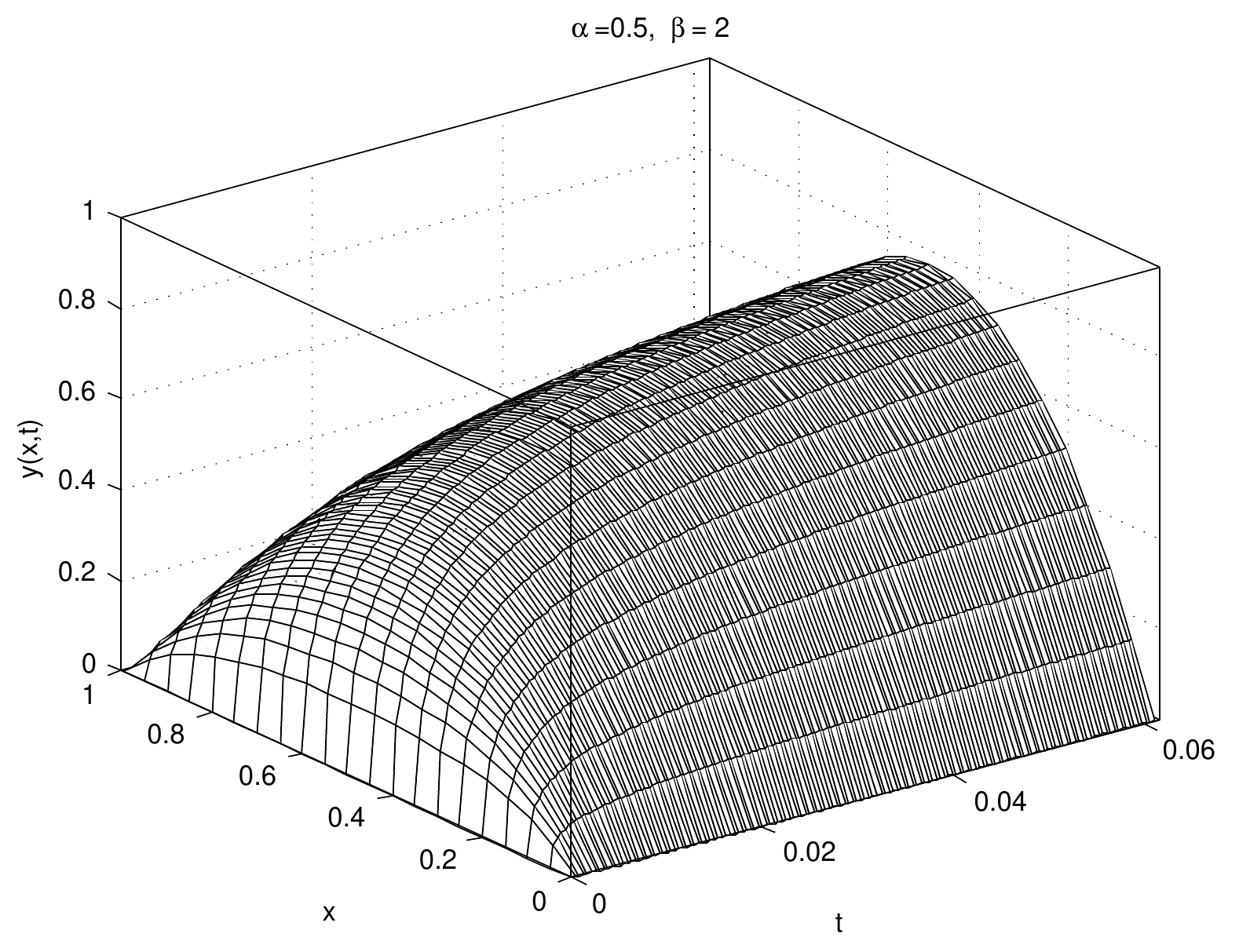} \hfill  \includegraphics[width=6cm]{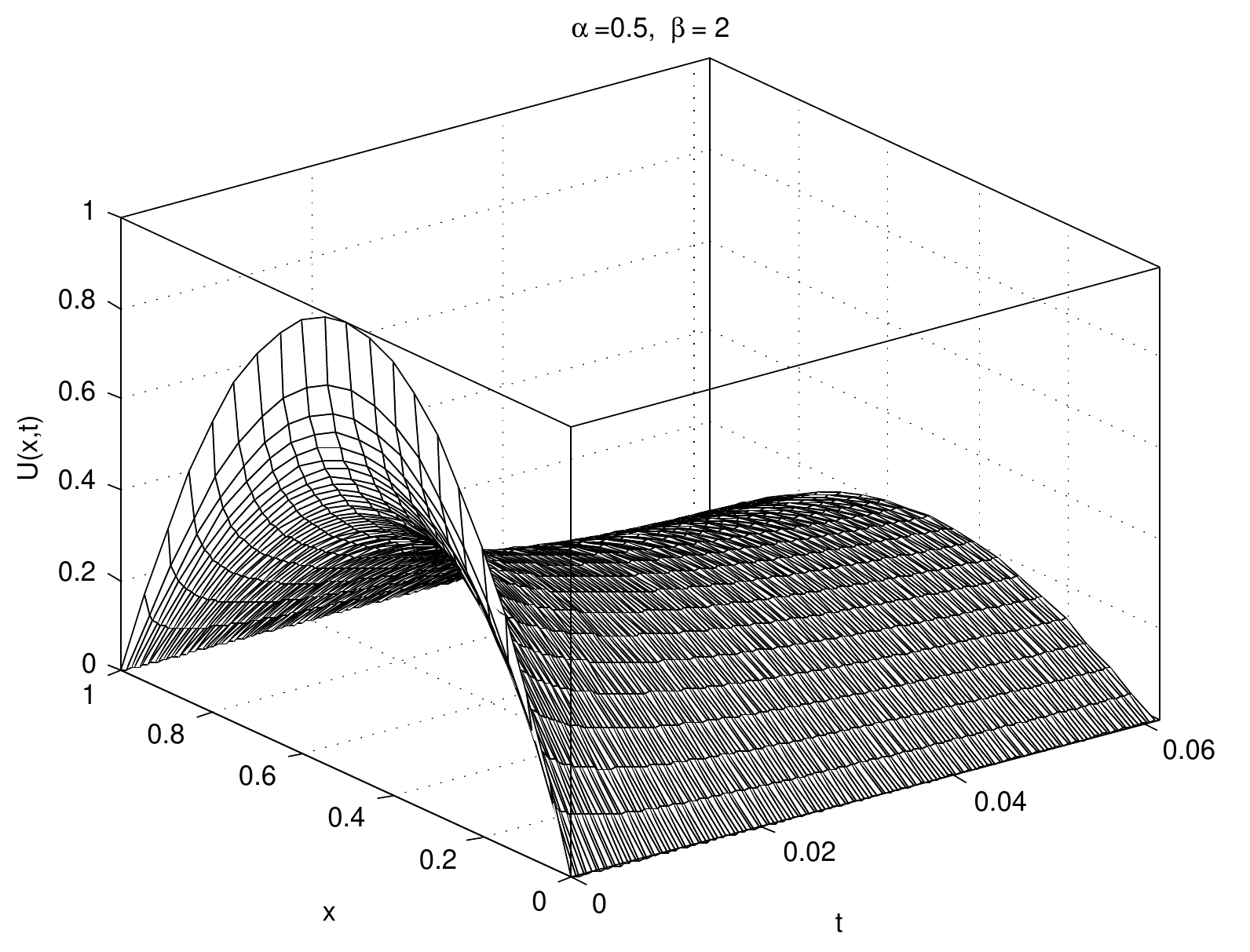}\\ 
\caption{Solutions $y(x,t)$ (left column) and $u(x,t)$ (right column) of Example 2, for $\alpha =1$ (top), $\alpha=0.7$ (middle) and $\alpha=0.5$ (bottom), with spatial step $h=0.05$ and time step $\tau=h^2/6$.}
\label{fig:example-2}
\end{center}
\end{figure}

The results of computations of $y(x,t)$ and then $u(x,t)$  for $\alpha =1$,  $\alpha = 0.7$, $\alpha = 0.5$ with $h=0.05$ and $\tau=h^2/6$  are shown in Fig.~ \ref{fig:example-2}.  
The structure of the matrix is the same as shown in Fig.~\ref{fig:matrix-structure-zoomed}.

Obviously, for $\alpha = 1$ we have the classical case and the same plots 
as in Fig.~\ref{fig:example-1}, and therefore Example 1 is a particular case of Example 2. 
As $\alpha$ goes to zero, the function $y(x,t)$ slowly tends to $u(x,0) = 4x(1-x)$ for all $t$.  
This is also not a surprize, because, indeed, for $\alpha=0$ the function $y(x,t)$ does not
depend on $t$ and therefore must satisfy 

\begin{displaymath}
y''(x) + 8 = 0, 
\quad
y(0) = y(1) = 0,
\end{displaymath} 

\noindent
which has the solution $y(x) = 4x(1-x)$.

It should be noted that almost the same problem as 
(\ref{eq:example-2-equation})--(\ref{eq:example-2-ic})
was numerically solved in \cite{Scherer:2008kk} using a very different approach.
The initial condition in  \cite{Scherer:2008kk} was  $u(x,0) = x(1-x)$.
Scaling the plots   in figures~1 and 2 in \cite{Scherer:2008kk} by the factor of 4,
we obtain the plots which are practically identical with our results for $u(x,t)$ 
shown in Fig.~\ref{fig:example-2}.
For this comparison we considered the shorter interval $0\leq t \leq 0.02$ 
used in \cite{Scherer:2008kk}.

\bigskip
\bigskip

\subsection{Example 3: Diffusion equation with spatial fractional derivative}

Let us now focus on the role of \emph{spatial} fractional derivative. 
For clarity, let us directly write the following analog of the problem 
(\ref{eq:example-1-equation-y})--(\ref{eq:example-1-ibc-y}) 
for determining the function $y(x,t)$:

\begin{equation} \label{eq:example-3-equation-y}
\frac{\partial y}{\partial t} -
\frac{\partial^\beta y}{\partial |x|^\beta} = f(x,t),  \qquad (\mbox{with } f(x,t) \equiv  8)
\end{equation}

\begin{equation} \label{eq:example-3-ibc-y}
y(0,t) = 0, \quad y(1,t) = 0; 
\qquad  \quad
y(x,0) = 0.
\end{equation}

\noindent
where $1 < \beta \leq 2$. 
The right-hand side is the same as in (\ref{eq:example-1-equation-y}), but instead of second order spatial derivative we deal with the Riesz-Caputo fractional derivative. 

The problem (\ref{eq:example-3-equation-y})--(\ref{eq:example-3-ibc-y}) can be discretized using the described method (see Fig.~\ref{fig:discretization-equation}), which gives

\begin{equation}\label{eq:example-3-discretization}
\Bigl\{
B_{n}^{(1)} \otimes E_m - \,
E_n \otimes R_{m}^{(\beta)}
\Bigr\} y_{nm} = f_{nm}
\end{equation}

\noindent
where $m$ is the number of spatial discretization intervals and $n$ is the number of time steps, 
and the corresponding rows and columns in the system (\ref{eq:example-3-discretization}) are removed with the help of eliminators. 

\begin{figure}[t]
\begin{center}
\includegraphics[width=6cm]{images/y-10-20} \hfill  \includegraphics[width=6cm]{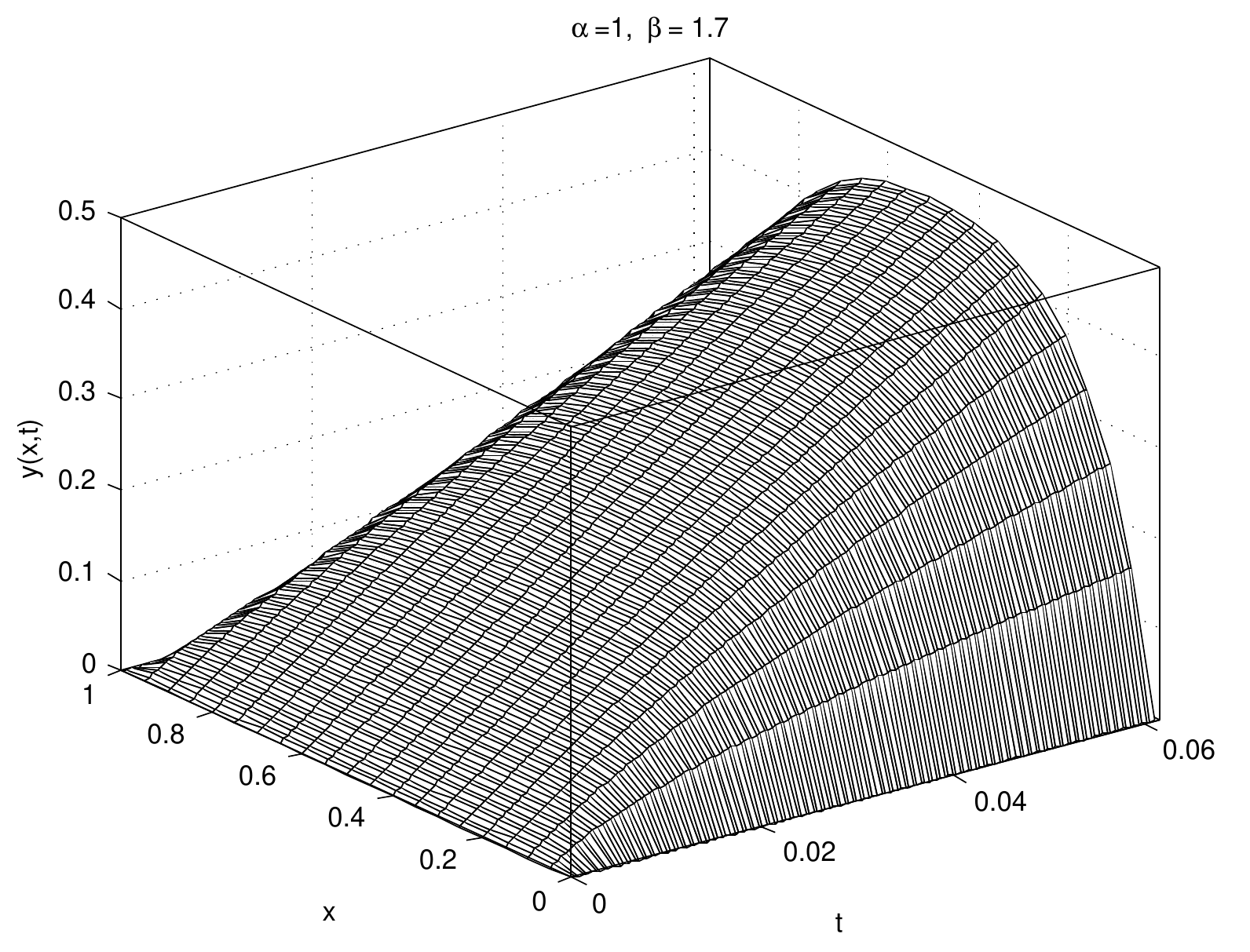}\\ 
\includegraphics[width=6cm]{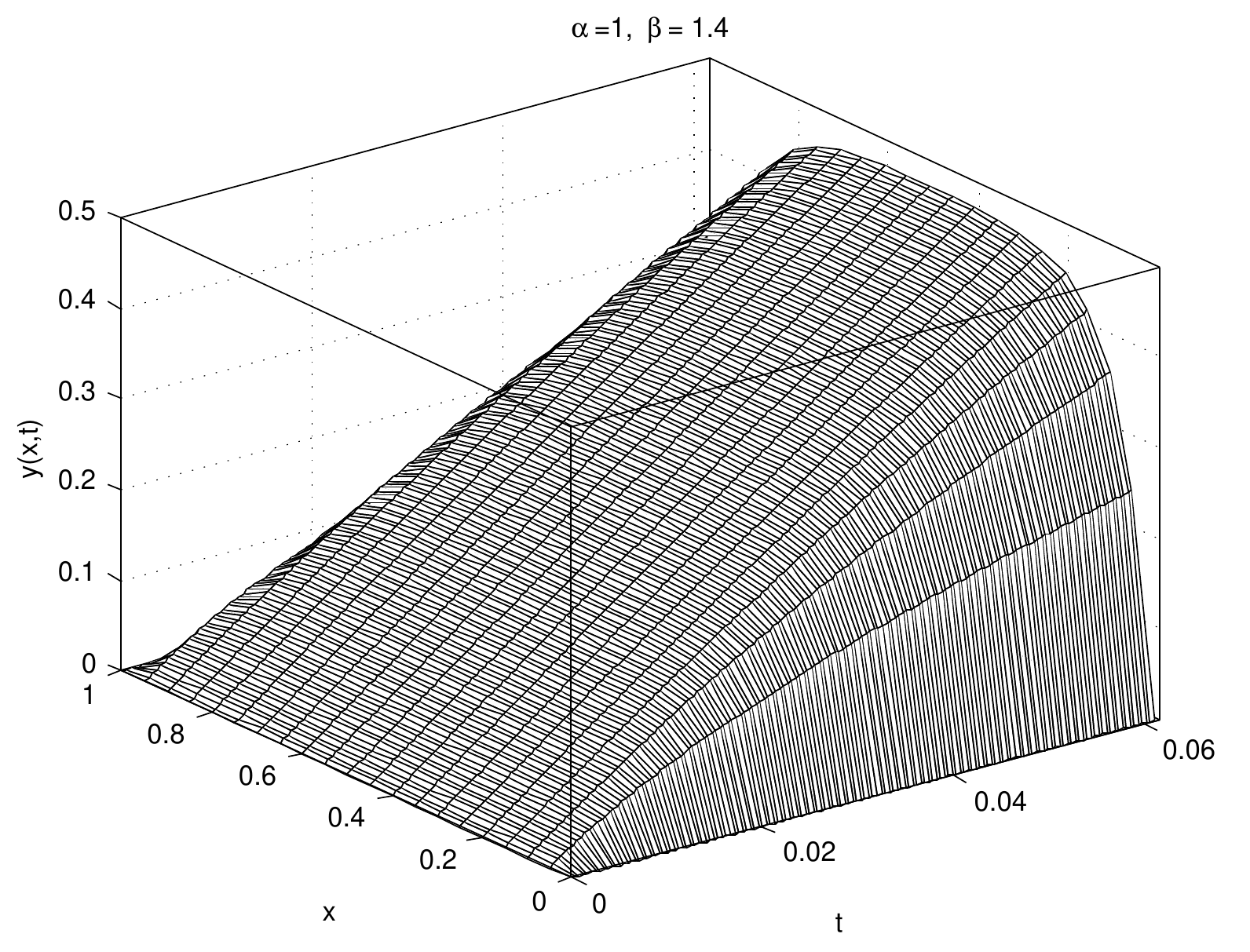} \hfill  \includegraphics[width=6cm]{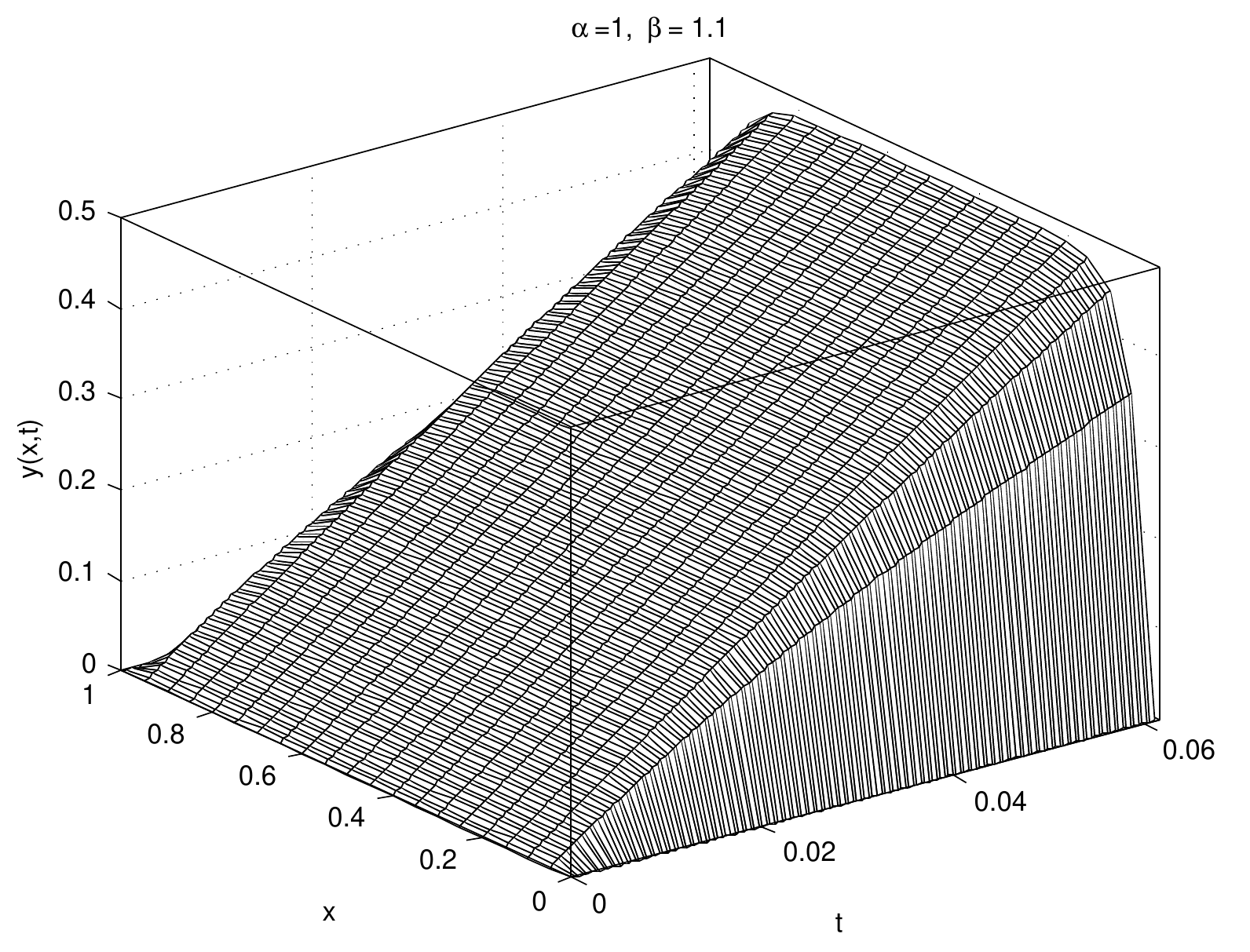}\\ 
\caption{Solutions $y(x,t)$ (left column) of Example 3, for $\beta =2$ (top-left), $\beta=1.7$ (top-right), $\beta=1.4$ (bottom-left), and $\beta=1.1$ (bottom-right),
with spatial step $h=0.05$ and time step $\tau=h^2/6$.}
\label{fig:example-3}
\end{center}
\end{figure}

The results of computations for four different values of $\beta$ are shown in Fig.~\ref{fig:example-3}.

\bigskip
\bigskip

\subsection{Example 4: General fractional diffusion equation}

Now we can illustrate that the method works also in the case when both derivatives
are of fractional order. Let us consider the most general situation:

\begin{equation} \label{eq:example-4-equation-y}
 _{0}^{C}\!D_{t}^{\alpha}y -
\frac{\partial^\beta y}{\partial |x|^\beta} = f(x,t),  \qquad (\mbox{with } f(x,t) \equiv  8)
\end{equation}

\begin{equation} \label{eq:example-4-ibc-y}
y(0,t) = 0, \quad y(1,t) = 0; 
\qquad  \quad
y(x,0) = 0.
\end{equation}

The right-hand side is the same as in (\ref{eq:example-1-equation-y}) and (\ref{eq:example-3-equation-y}), but now \emph{both} derivatives are allowed to be of non-integer order.

The problem (\ref{eq:example-4-equation-y})--(\ref{eq:example-4-ibc-y}) can be discretized using the described method (see Fig.~\ref{fig:discretization-equation}), which gives

\begin{equation}\label{eq:example-4-discretization}
\Bigl\{
B_{n}^{(\alpha)} \otimes E_m - \,
E_n \otimes R_{m}^{(\beta)}
\Bigr\} y_{nm} = f_{nm}
\end{equation}

\noindent
where $m$ is the number of spatial discretization intervals and $n$ is the number of time steps, 
and the corresponding rows and columns in the system (\ref{eq:example-4-discretization}) are as in all previous examples removed with the help of eliminators.

\begin{figure}
\begin{center}
\includegraphics[width=4cm]{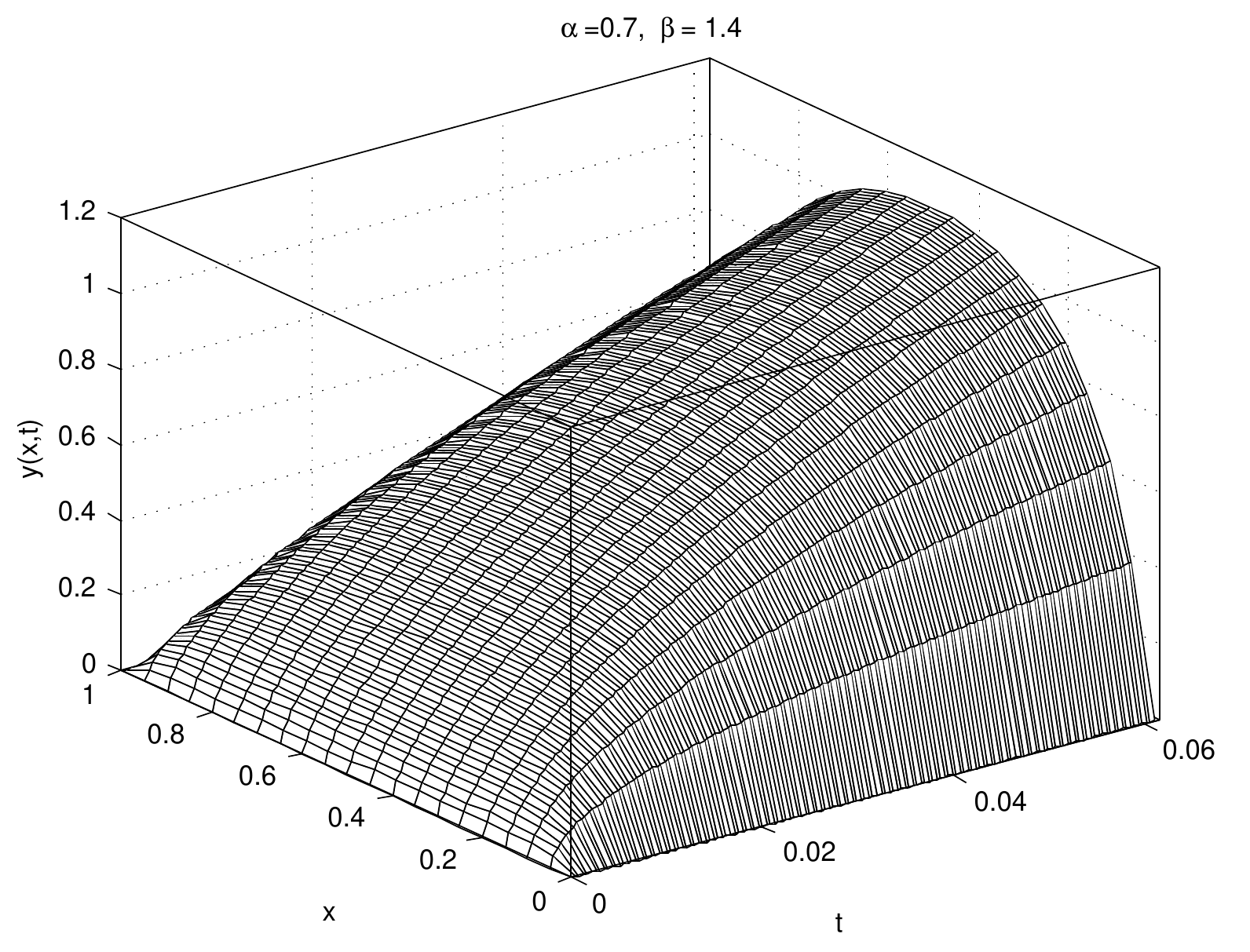} \hfill  
\includegraphics[width=4cm]{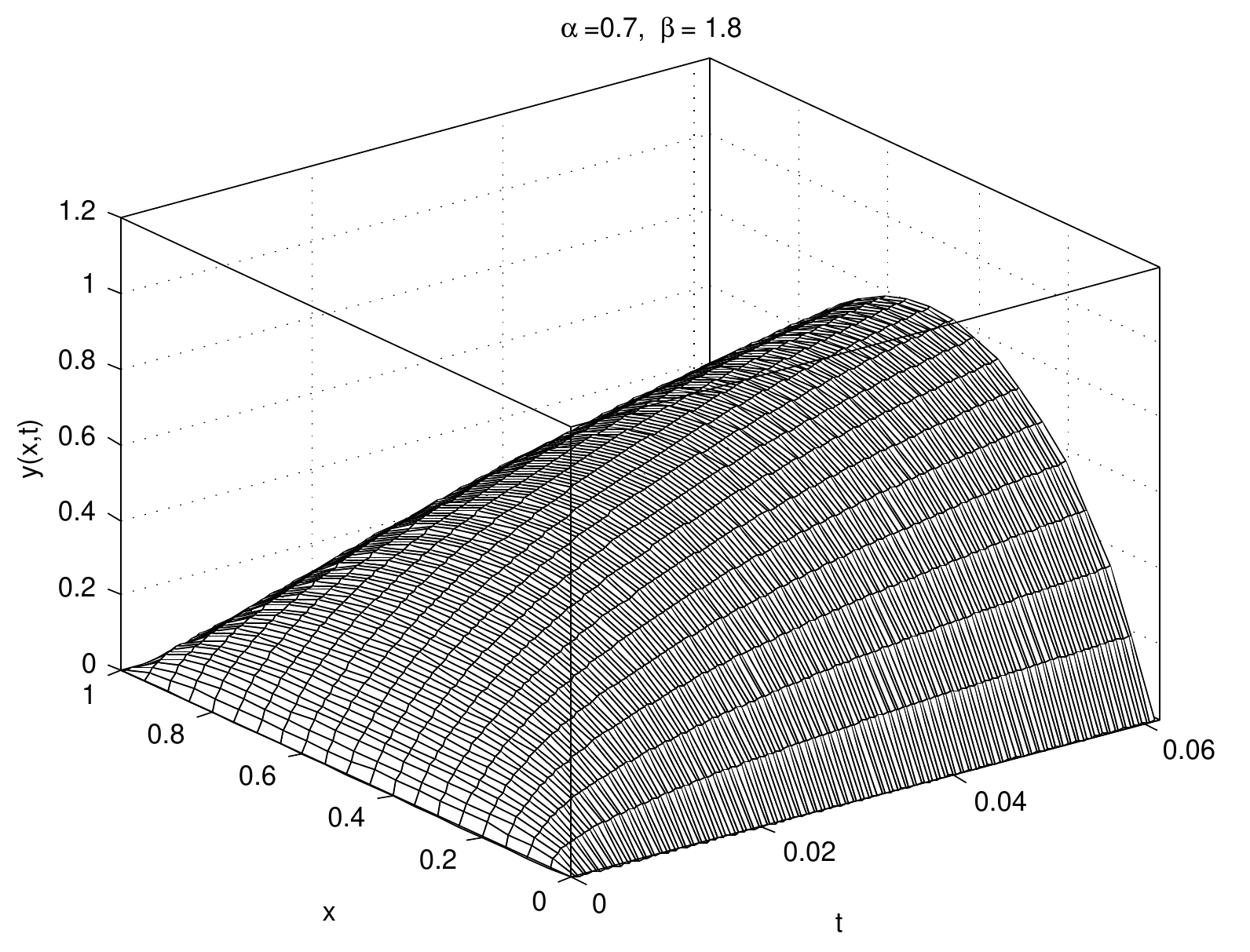} \hfill
\includegraphics[width=4cm]{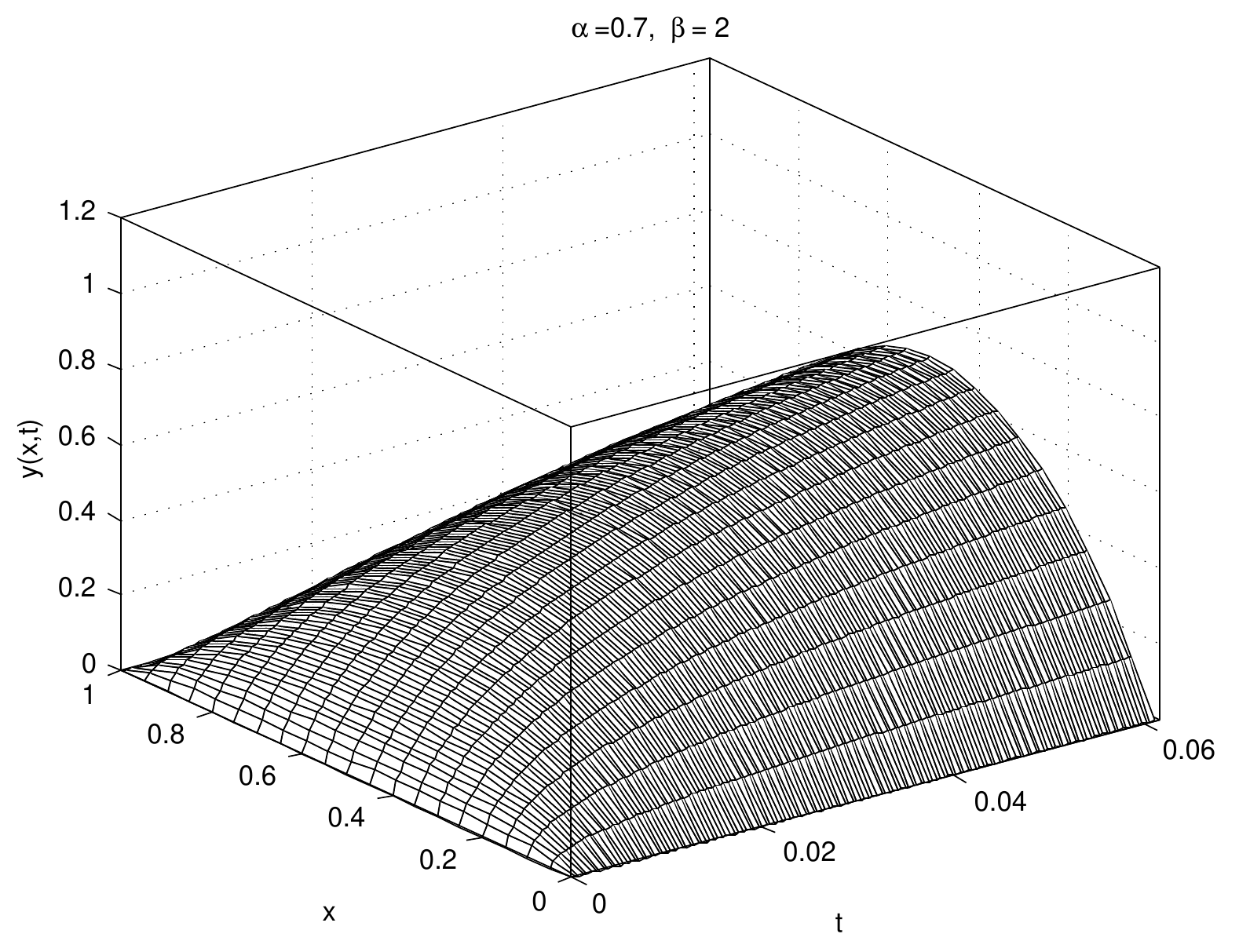} \hfill  
\caption{Solutions $y(x,t)$ of Example 4, for $\alpha=0.7$ and $\beta =1.4$ (left), 
$\alpha=0.7$ and $\beta=1.8$ (middle), $\alpha=0.7$ $\beta=2$ (right),
with spatial step $h=0.05$ and time step $\tau=h^2/6$.}
\label{fig:example-4}
\end{center}
\end{figure}

The results of computations for some sample combinations of non-integer orders $\alpha$ and different values of $\beta$ are shown in Fig.~\ref{fig:example-4}.

\bigskip
\bigskip

\subsection{Example 5: Fractional diffusion equation with delay}

Finally, let us consider the equation with two Caputo fractional-order time derivatives,
of which one is with delay $\delta$ (we do not go into the physical interpretation of this equation,
because physical interpretation of a delayed fractional derivative is not known so far,
but use it for demonstrating how broad can be the field of application of our approach):

\begin{equation}\label{eq:example-5-equation-y}
\frac{1}{2}
\left\{
\; _{0}^{C}\!D_{t}^{\alpha}y + \; 
 _{0}^{C}\!D_{t-\delta}^{\gamma}y
\right\}
 - \frac{\partial^\beta y}{\partial |x|^\beta}  = f(x,t)
  \qquad (\mbox{with } f(x,t) \equiv  8)
\end{equation}

\begin{equation} \label{eq:example-5-ibc-y}
y(0,t) = 0, \quad y(1,t) = 0 \quad  y(x,0) = 0
\end{equation}

Obviously, for $\gamma = \alpha$ and $\delta=0$ we have the equation considered in Example 2. 
Let us select the discretization step so that $\delta$ is a multiple of the time step $\tau$: 
$\delta = k \tau$. 
Then the problem (\ref{eq:example-5-equation-y})--(\ref{eq:example-5-ibc-y}) can be discretized using the described method (see Fig.~\ref{fig:discretization-equation} and the equation (\ref{eq:shift-ne})), which gives:

\begin{equation}\label{eq:example-5-discretization}
\Bigl\{
\frac{1}{2}
\left(
 B_{n}^{(\alpha)} \otimes E_m + \,
 \,
_{+k}B_{n}^{(\gamma)} \otimes E_m 
\right)
- \,
\, E_n \otimes R_{m}^{(\beta)}
\Bigr\} y_{nm} = f_{nm}
\end{equation}

\begin{displaymath}
_{+k}B_{n}^{(\gamma)}  = S_{n+1,\dots, n+k} \, E_{n+k,k}^{+} \, B_{n+k}^{(\gamma)}  \, E_{n+k,k}^{+} \, S_{1,\ldots, k}^T 
\end{displaymath}

\noindent
where, as above,
 $m$ is the number of spatial discretization intervals and $n$ is the number of time steps, 
$k$ is the number of time steps corresponding to the delay $\delta$,
and the appropriate rows and columns in the system (\ref{eq:example-5-discretization}) are as in all previous examples are to removed with the help of eliminators.

\begin{figure}[ht]
\begin{center}
\includegraphics[width=6cm]{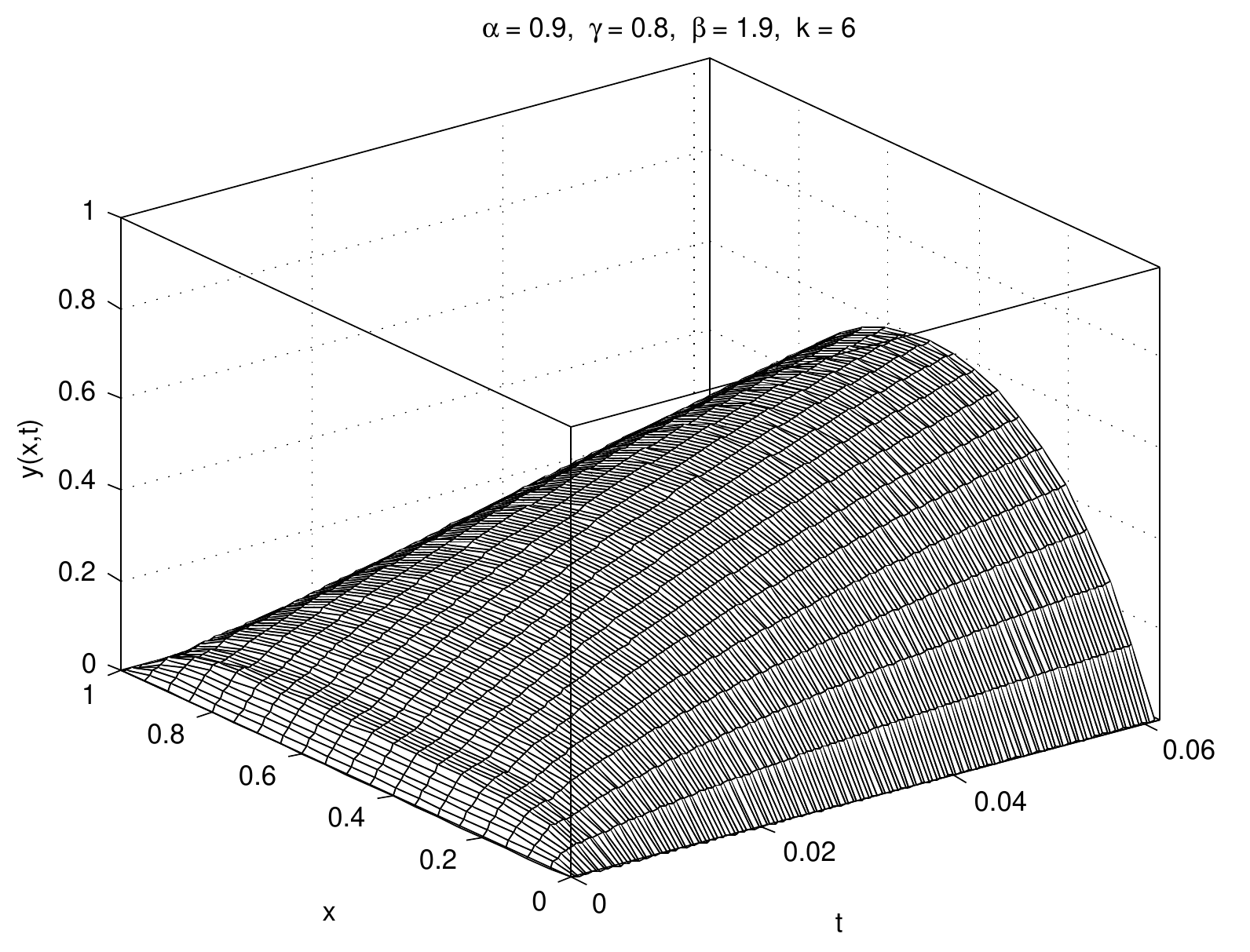} \hfill  
\includegraphics[width=6cm]{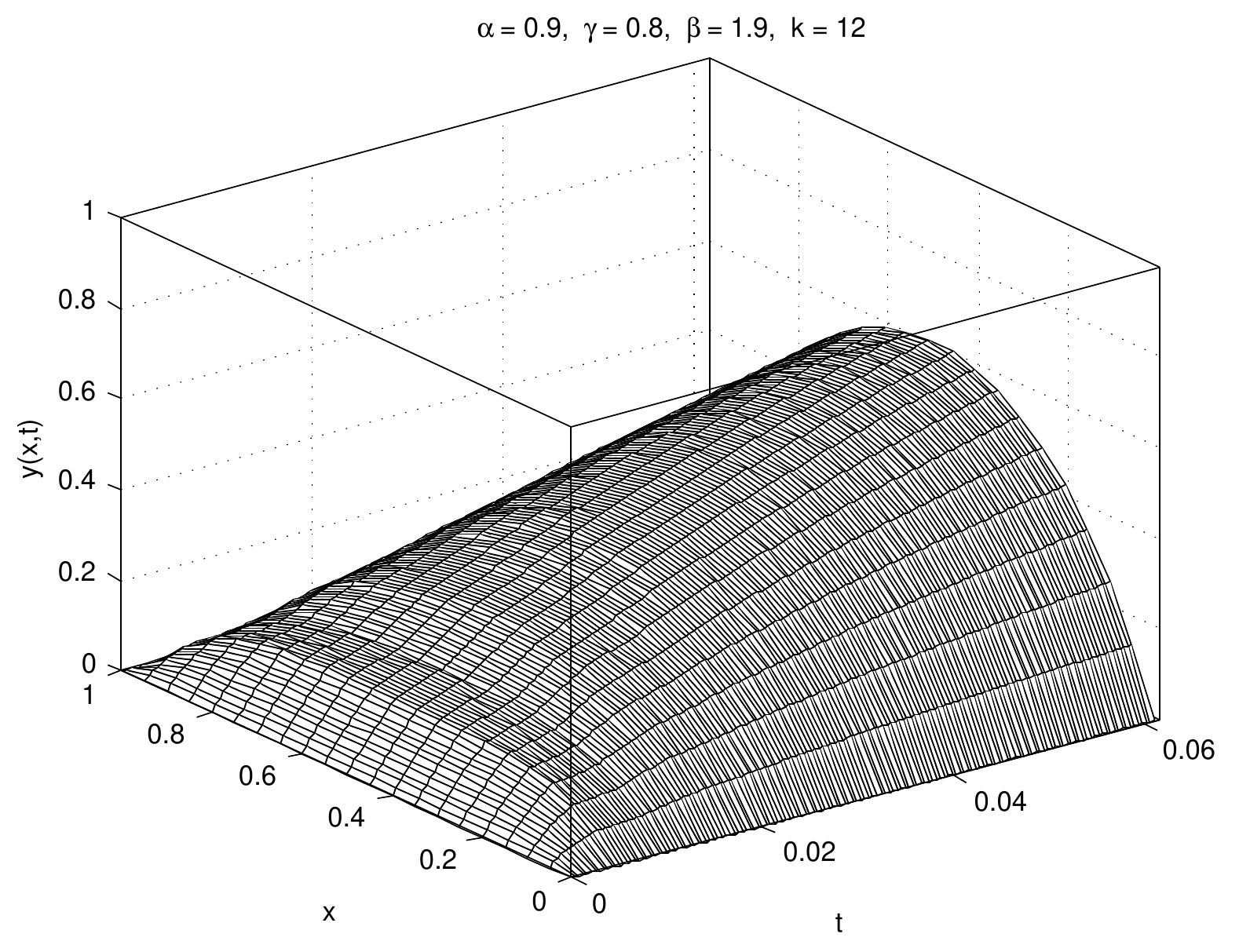} \\
\includegraphics[width=6cm]{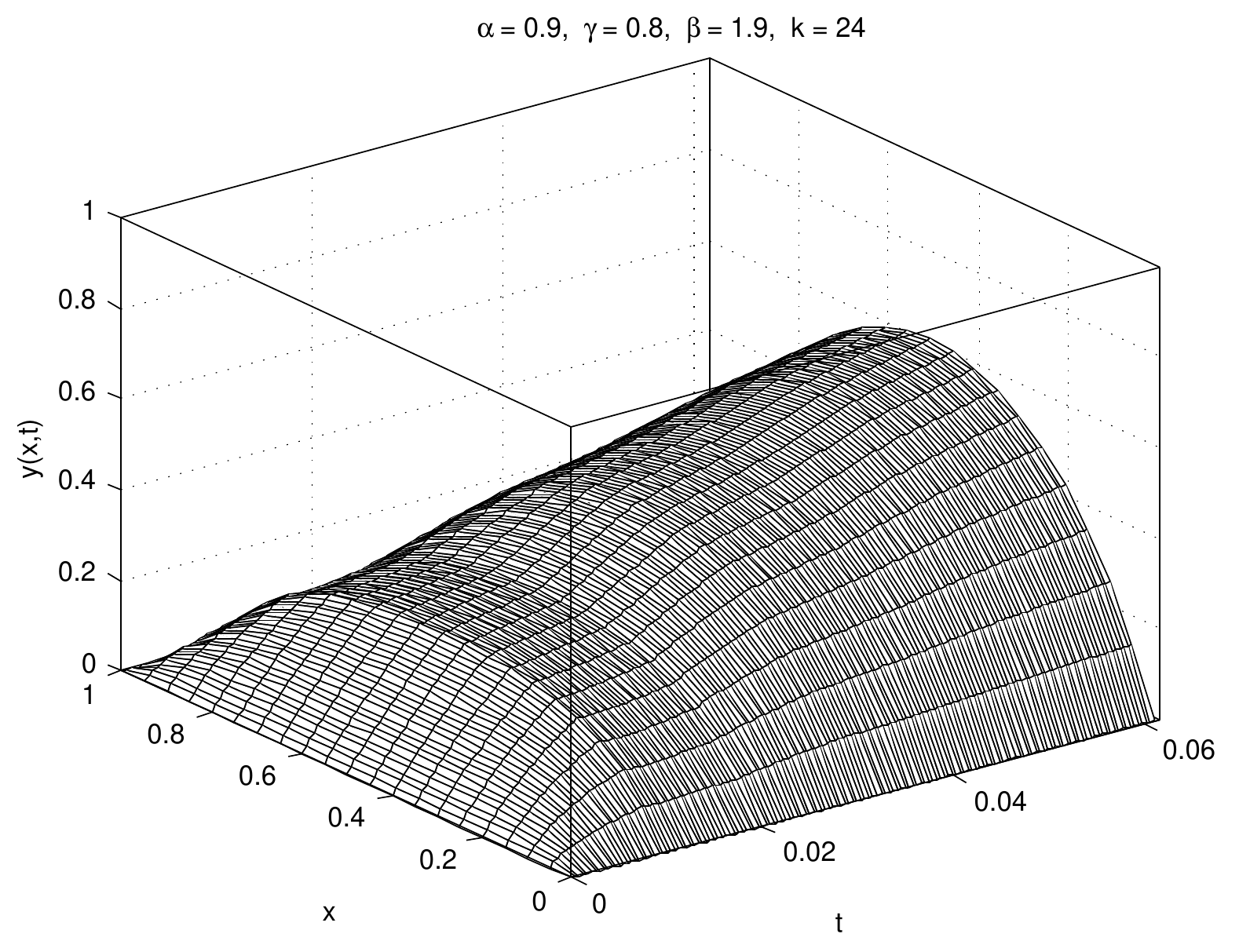} \hfill
\includegraphics[width=6cm]{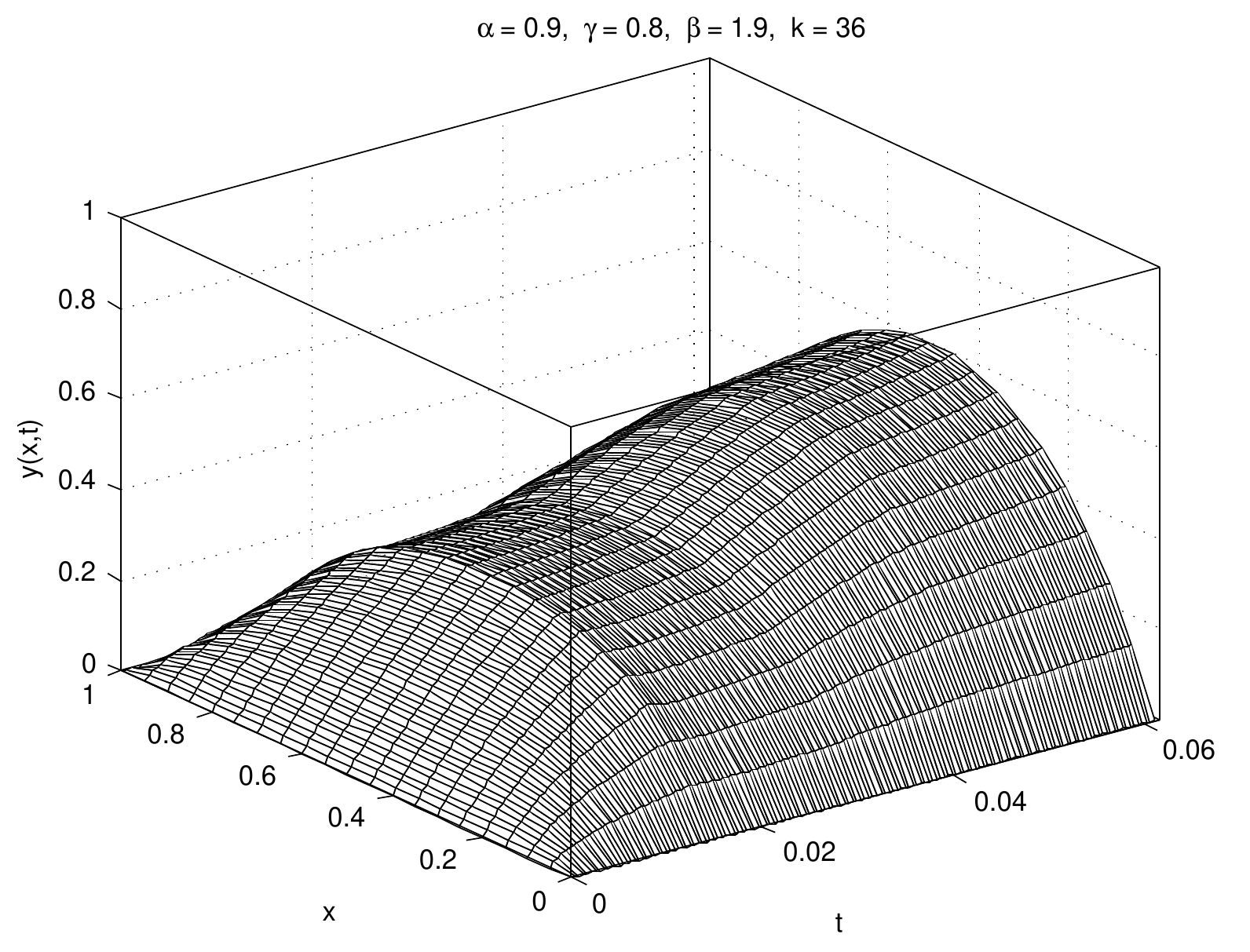} 
  \caption{Solutions $y(x,t)$ (left column) of Example 5, for $\alpha=0.9$, $\gamma$ = 0.8, $\beta =1.9$, for delays $\delta_k = k \tau$, k = 6, 12, 24, 36.   }
\label{fig:example-5}
\end{center}
\end{figure}

The results of computations for a sample combination of non-integer orders $\alpha$, $\beta$ and $\gamma$ and some delays $\delta$  represented by the parameter $k$ are shown in Fig.~\ref{fig:example-5}.

\section{Conclusion}

The suggested method represents a unifying approach to numerical solution of partial differential equations of both integer and non-integer order, including equations with delays.

For the sake of clarity, in this article we considered the case of one spatial variable.
However, the suggested method can be easily extended to the case of two and three spatial variables by repeatedly applying the triangular strip matrix representations of fractional-order operators in combination with the Kronecker matrix product. 

The problems considered in this article are linear, so the resulting systems of algebraic equations are linear as well. However, the suggested approach can be extended to the case of nonlinear problems, too. 

The suggested method can be used also for solving partial fractional differential equations appearing from the Laplace equation by replacing second order spatial derivatives with fractional Riesz derivatives. 

The suggested method can be used also for partial fractional FDEs of variable and distributed order(s) and for equations with delays.

\section*{Acknowledgments}

This work could not be done without academic exchanges supported by the National Scholarship Program of the Slovak Republic (visit of A.~Chechkin to Kosice, Slovakia), Summer Fellowships Program of the Department of Electrical and Computer Engineering of the Utah State University (visit of I.~Podlubny to Logan, USA), and the support grant SAB2006-0172 from the Ministry of Education of Spain (visit of I.~Podlubny to the University of Extremadura, Badajoz, Spain). 

The authors are grateful to Professor Francesco Mainardi and to anonymous referees for their valuable comments and suggestions. 

\section*{Appendix: sample evaluation of the symmetric Riesz fractional derivative}

For $\phi (x) = x (1-x)$ and the order of differentiation $1<\beta< 2$ 
the left-sided Riemann-Liouville fractional derivative (\ref{eq:RL-definition-left})
 of the function $\phi (x)$ is

\begin{equation}
    	_{0}D_{x}^{\beta} \phi(x) 
	=
	\frac{x^{1-\beta}}{\Gamma(2-\beta)}
	- 
	\frac{2 \, x^{2-\beta}}{\Gamma(3-\beta)}.   
\end{equation}

Similarly, the right-sided Riemann-Liouville derivative (\ref{eq:RL-definition-right}) of $\phi(x)$ is

\begin{equation}
    	_{x}D_{1}^{\beta} \phi(x) 
	=
	\frac{(1-x)^{1-\beta}}{\Gamma(2-\beta)}
	- 
	\frac{2 \, (1-x)^{2-\beta}}{\Gamma(3-\beta)}.    
\end{equation}

Therefore, the symmetric Riesz fractional derivative (\ref{eq:Riesz-derivative-definition}) of the function $\phi(x)$ is:

\begin{eqnarray}
	\frac{d^\beta \phi}{d |x|^\beta} 
	 & = &
	\frac{1}{2}
	\left\{ \, 
		_{0}D_{x}^{\beta} \phi(x) 
		+ \,
		_{x}D_{1}^{\beta} \phi(x) 
	\right\} \\ 
	& = & 
	\frac{x^{1-\beta} + (1-x)^\beta}
		{2  \, \Gamma(2-\beta)}
		-
		\frac{x^{2-\beta} + (1-x)^{2-\beta}}
			{ \Gamma(3-\beta)}.
\end{eqnarray}

\bibliography{nummetfracdiff-2a}

\end{document}